%% file: arXiv_qHP_MFIE_Mueller.tex
\documentclass[journal]{IEEEtran}

\usepackage{cite}
\usepackage{amssymb,amsmath, amsfonts, amsthm}
\usepackage[caption=false,font=footnotesize]{subfig}
\usepackage{tikz,pgfplots}
\usetikzlibrary{shapes.geometric, decorations.pathmorphing, positioning, arrows, arrows.meta, patterns, shapes, calc, fit, backgrounds, external}

\usepackage{graphicx}
\graphicspath{{./figures/}}
\DeclareGraphicsExtensions{.pdf,.jpeg,.png,.eps}
\usepackage[dvipsnames]{xcolor}

\input{lvc_style}

\begin{document}

\title{Suppressing Low-Frequency Cancellation Errors in Far-Field Computation with the MFIE, M\"{u}ller, and Multitrace M\"{u}ller Integral Equations}

\author{Chien~V.~Le,~\IEEEmembership{Member,~IEEE,} and~Kristof~Cools,~\IEEEmembership{Member,~IEEE}
\thanks{Manuscript received April 19, 2005; revised August 26, 2008; accepted 19 May 2014. Date of publication 12 June 2014; date of current version 9 July 2014. This work was supported by the European Research Council (ERC) under the European Union's Horizon 2020 Research and Innovation Program under Grant 101001847. \textit{(Corresponding author: Chien V. Le.)}}%
\thanks{The authors are with the IDLab, Department of Information Technology at Ghent University -- imec, 9000 Ghent, Belgium (e-mail: vanchien.le@ugent.be, kristof.cools@ugent.be).}
}


\maketitle
 
\begin{abstract}
     This paper presents robust numerical formulations for suppressing low-frequency cancellation errors in the magnetic field integral equation (MFIE) for perfect electric conductors and the M\"uller integral equations for homogeneous and composite dielectric objects. The proposed approach combines quasi-Helmholtz projectors, which separate the Helmholtz components of the surface currents, with a tailored rescaling procedure that preserves the components essential for accurate far-field evaluation from finite-precision cancellation. The resulting formulations retain the conditioning properties of the mixed MFIE and M\"uller discretizations in the dense-mesh and low-frequency regimes. Numerical experiments involving simply and multiply connected scatterers, composite dielectric structures, multiscale and high-contrast configurations demonstrate that the proposed methods yield accurate far-field predictions across a broad frequency range, from the full-wave regime to extremely low frequencies.
\end{abstract}



\section{Introduction}

\IEEEPARstart{B}{oundary} integral equations (BIEs) are widely used in the numerical modeling of electromagnetic scattering by conducting and dielectric bodies. Classical formulations include the electric field integral equation (EFIE) and the magnetic field integral equation (MFIE) for perfect electric conductors (PECs), together with the Poggio--Miller--Chang--Harrington--Wu--Tsai (PMCHWT) and M\"uller equations for dielectrics. These current-based formulations recast the scattering problems in terms of equivalent surface current densities, which represent the tangential traces of the electromagnetic fields on the boundary of the scatterer.

A key analytical tool for characterizing the behavior of surface currents is their Helmholtz decomposition into irrotational (non-solenoidal or star), harmonic (global-loop), and non-harmonic solenoidal (local-loop) components
\begin{equation}
    \label{eq:decomposition}
    \jb = \jb^{\Sm} + \jb^{\textrm{H}} + \jb^{\Lambda},
\end{equation}
respectively. In the low-frequency regime, i.e., as the angular frequency $\omega \to 0$, these Helmholtz components may exhibit different asymptotic scalings. Their scalings depend on the problem under consideration, namely PEC or dielectric scattering, the excitation, and the scatterer topology, which determines whether global-loop currents can occur. They are however independent of the formulation employed. As an example, Table~\ref{tab:scaling}, adapted from \cite{HEA+2023}, summarizes the low-frequency scalings of the Helmholtz components of the electric surface current density induced on a PEC scatterer by different excitations. Without loss of generality, the magnitudes of the components can be expressed as
\[
     \abs{\jb^{\Sm}} = \OO\paren{\omega^\alpha}, \qq 
     \abs{\jb^{\textrm{H}}} = \OO\paren{\omega^\beta}, \qq 
     \abs{\jb^{\Lambda}} = \OO\paren{\omega^\gamma}, 
\]
where $\alpha, \beta, \gm \in \Z$. A vanishing component is represented by the convention $0 = \OO(\omega^\infty)$. As a result, the relative contributions of the components may differ by several orders of magnitude in the low-frequency regime.

\begin{table*}[!t]
    \renewcommand{\arraystretch}{1.5}
    \caption{Low-frequency scaling of the Helmholtz components of the electric surface current in PEC scattering, as powers of $\omega$.}
    \label{tab:scaling}
    \centering
    \begin{tabular}{| c c  c | c c c c c c c|}
      \hline
      \multicolumn{3}{|c|}{Surface Current $\jb$} & \multicolumn{7}{|c|}{Excitations} \\
      \hline
      Component & & Exponent & Plane wave & Hertzian dipole & elec. ring current & mag. dipole & mag. ring current & ind. gap & cap. gap \\ \hline
      Non-solenoidal & $\jb^{\Sm}$ & $\alpha$ & 1 & 0 & 2 & 1 & 1 & 1 & 1 \\
      Harmonic & $\jb^{\textrm{H}}$ & $\beta$ & 0 & 0 & 0 & -1 & 1 & -1 & 1 \\
      Local-loop & $\jb^{\Lambda} $ & $\gamma$ & 0 & 0 & 0 & -1 & 1 & -1 & 1 \\ \hline
    \end{tabular}
\end{table*}

This imbalance gives rise to a fundamental numerical issue when solving BIEs at low frequencies in finite-precision arithmetic. When all components are stored together within a single floating-point variable, those with larger exponents (i.e., smaller magnitudes) are susceptible to loss of significant digits against dominant components. Although such errors may remain less pronounced in the surface currents and reconstructed near fields, they become critical when evaluating the far field.

Indeed, the far-field pattern associated with a surface current density $\jb$ on the boundary $\Gm$ is stably computed as follows \cite{ADC+2021}:
\begin{equation}
    \label{eq:farfield}
    \begin{aligned}
        \Eb^{\text{far}}(\wh{\xb}) & = \int_\Gm e^{-\iota \kappa_0 \, \wh{\xb} \cdot \yb} \, \jb^{\Sm}(\yb) \ds_{\yb} \\
        & \, + \int_\Gm \paren{e^{-\iota \kappa_0 \, \wh{\xb} \cdot \yb} - 1} \paren{\jb^{\textrm{H}}(\yb) + \jb^{\Lambda}(\yb)} \ds_{\yb},
    \end{aligned}
\end{equation}
where $\iota = \sqrt{-1}$ is the imaginary unit, $\kappa_0 = \omega/c_0$ and $c_0 = 1/\sqrt{\mu_0 \epsilon_0}$ are respectively the wavenumber and the wave speed in the background medium with permittivity $\epsilon_0$ and permeability $\mu_0$, and $\wh{\xb} = \xb / \abs{\xb}$ is the direction.
Since
\[
    e^{-\iota \kappa_0 \, \wh{\xb} \cdot \yb} = \OO(1), \qqq e^{-\iota \kappa_0 \, \wh{\xb} \cdot \yb} - 1 = \OO(\omega),
\]
as $\omega \to 0$, the non-solenoidal component $\jb^{\Sm}$ may contribute non-negligibly to the far field even when it is asymptotically smaller than the solenoidal components. In particular, $\jb^{\Sm}$ contributes at the same asymptotic order as $\jb^{\textrm{H}} + \jb^{\Lambda}$ when
\begin{equation}
    \label{eq:albtgm}
    \alpha = \min(\beta, \gm) + 1.
\end{equation}
However, this is precisely the regime in which $\jb^{\Sm}$ is susceptible to loss of significant digits. As a result, the far-field evaluation becomes highly unreliable, often leading to large inaccuracies or even incorrect radiation patterns.

This phenomenon, commonly referred to as low-frequency cancellation or low-frequency catastrophe, is well known and has been extensively studied in literature \cite{ZCC+2003,VGG+2013,LSC2018}, particularly in the context of the EFIE \cite{QC2010,ACB+2013,HEA+2023b}, combined field integral equations (CFIEs) \cite{MBC+2020,LCA+2024}, and PMCHWT formulations \cite{GAM+2017,BMC+2017,LMA+2026}. It stems from the interplay between the different asymptotic scalings of the surface current components and the limitations of finite-precision arithmetic. Several remedies have been proposed, including potential-based formulations \cite{VFG+2016,LSC2018}, current-charge formulations \cite{QC2010,CAY+2015}, and quasi-Helmholtz decompositions combined with suitable component-wise rescaling schemes \cite{ACB+2013,MBC+2020,ADC+2021}. 

This work addresses low-frequency cancellation errors in the MFIE and the M\"uller equations. These formulations are particularly attractive targets for three reasons. First, as second-kind integral equations, they are inherently well-conditioned in the dense-mesh and low-frequency regimes. For the MFIE, this favorable conditioning is restricted to simply connected geometries and to wavenumbers sufficiently far from spurious resonances \cite{CAO+2009a}. Consequently, in the absence of topology-induced breakdown and resonant instability, finite-precision cancellation becomes the principal obstacle to the reliable use of these formulations at low-frequencies. Despite its practical importance, this issue has received comparatively little attention. In \cite{VGG+2013}, low-frequency cancellation in the MFIE was mitigated by solving an additional boundary integral equation for the surface charge. The authors of \cite{MBC+2020} instead considered a symmetrized formulation obtained by multiplying the exterior MFIE by an interior counterpart. The resulting symmetric structure enables cancellation errors in the MFIE and the corresponding CFIE to be suppressed through quasi-Helmholtz decompositions combined with an appropriate rescaling procedure, without degrading the conditioning of the formulations. To the best of the authors' knowledge, no analogous method have been proposed for the M\"uller equation.

Second, the MFIE and M\"uller formulations have historically been considered less accurate than their first-kind counterparts, such as the EFIE and PMCHWT formulations, particularly for far-field computations \cite{ZCC+2003,YTJ2005,YTJ2008,UTR2011,KE2023}. However, the development of mixed Petrov--Galerkin discretization schemes employing dual test functions, such as the Buffa--Christiansen (BC) basis functions \cite{BC2007}, has significantly improved their accuracy, in many cases to a level comparable to that of first-kind equations \cite{CAD+2011,YJN2011,YJN2013,YKJ2016,BCA+2014}. In addition, these mixed discretizations correctly reproduce the low-frequency scaling of the continuous solution at the discrete level \cite{BCA+2014}, which is a prerequisite for the successful recovery of current components affected by cancellation.

Third, although first-kind formulations have been extensively stabilized against low-frequency cancellation errors, low-frequency and dense-mesh breakdown, and late-time instability, the resulting remedies often entail a significant increase in computational cost. This additional cost is mainly due to the larger number of operator matrices that must be assembled and applied, and becomes especially relevant for composite configurations and time-domain simulations. By focusing on the MFIE and M\"uller equations, we aim to provide a robust and comparatively simple alternative for electromagnetic scattering spanning from moderately high to extremely low frequencies, enabling accurate far-field evaluation without introducing unnecessary additional complexity.

In this work, we employ quasi-Helmholtz projectors to separate the Helmholtz components of the surface currents and a tailored rescaling strategy to protect the relevant components from loss of significant digits. This approach avoids the introduction of auxiliary charge equations as in \cite{VGG+2013}, which increases the number of unknowns and complicates the resulting system. It also avoids the construction of an additional complementary operator to symmetrize the equation as in \cite{MBC+2020}, which likewise leads to a more complicated formulation. Moreover, for the M\"{u}ller equation, it is not clear how such a complementary operator should be defined. Instead, we adapt the strategy proposed in \cite{BMC+2017} for the PMCHWT formulation. The key observation in \cite{BMC+2017} is that the standard left quasi-Helmholtz rescaling, which is effective for the EFIE in \cite{ACB+2013}, does not necessarily yield well-conditioned systems for other formulations in the low-frequency regime. In particular, when the geometry is multiply connected, applying the rescaling operator directly to the left of the discretized system may deteriorate the conditioning. To avoid this difficulty, the matrix system is first multiplied by the inverse Gram matrix. This maps the system back to the coefficient space, where the rescaling operator can be applied more appropriately and may recover the desired conditioning behavior. This rescaling strategy provides a simpler way to regularize the MFIE. On multiply connected surfaces or at wavenumbers near resonances, it does not render the MFIE fully well-conditioned, since the system remains affected by the inherent topology-induced breakdown and resonant instability \cite{CAO+2009a}. The resulting formulation also does not preserve the nullspace of the static MFIE operator, since this nullspace is highly perturbed by quadrature errors. Nevertheless, this strategy paves the way for the low-frequency regularization of the M\"{u}ller equations. The proposed method is particularly advantageous for composite configurations involving multiple objects and for time-domain simulations, where existing approaches are either not directly applicable or would lead to prohibitively complicated systems. 

As this work focuses on the accurate evaluation of the far field, we preserve the solution components that are essential to far-field computation but are susceptible to loss of significant digits. Specifically, we consider the regime in which the condition \eqref{eq:albtgm} is satisfied. Among the excitations listed in Table~\ref{tab:scaling}, this regime arises only for plane-wave excitation ($\alpha = 1$, $\beta = \gm = 0$). The other excitation types in Table~\ref{tab:scaling} either do not induce cancellation errors (i.e., $\alpha = \beta = \gamma$, such as the Hertzian dipole, magnetic ring current, and capacitive gap), or affect components that do not contribute significantly to the far field (i.e., $\alpha > \min(\beta, \gamma) + 1$, such as the electric ring current, magnetic dipole, and inductive gap). The formulations developed in this paper are therefore specialized to plane-wave excitation. For a more general excitation-aware and self-adaptive framework that recovers all solution components, including those relevant to other quantities of interest, we refer the reader to the exhaustive analysis in \cite{HEA+2023}.


The remainder of this paper is organized as follows. Section~\ref{sec:MFIE} introduces the proposed low-frequency regularization of the MFIE. Section~\ref{sec:Mueller} develops the corresponding regularization of the M\"uller equation for homogeneous dielectric objects. Section~\ref{sec:MT_Mueller} extends the proposed approach to general composite penetrable structures using the recently developed global multitrace (MT) M\"uller formulation in \cite{LC2026}. Implementation details and computational complexity are discussed in Section~\ref{sec:implementation}. Section~\ref{sec:results} presents numerical results that validate the proposed methods. Finally, conclusions are drawn in Section~\ref{sec:conclusion}.

\section{Magnetic Field Integral Equation}
\label{sec:MFIE}

\subsection{Formulation and Discretization}
\label{subsec:MFIE_formulation}

We consider a bounded and connected PEC $\Om \sst \R^3$ whose boundary $\Gamma$ is Lipschitz continuous and closed. The outward normal on $\Gm$ is denoted by $\nv$. The PEC is immersed in a homogeneous background $\Om_0 := \R^3 \setminus \ovl{\Om}$ with permittivity $\epsilon_0$ and permeability $\mu_0$. The scatterer $\Om$ is illuminated by a time-harmonic incident plane wave $(\eb^{inc}, \hb^{inc})$ with angular frequency $\om > 0$, originating in $\Om_0$ and satisfying the Maxwell equations there. The electric surface current density $\jb$ induced on $\Gm$ satisfies the MFIE
\begin{equation}
\label{eq:MFIE}
    \paren{\dfrac{1}{2} I - K^{(\kappa_0)}} \jb(\xb) = \nv \times \hb^{inc}(\xb).
\end{equation}
Here, $I$ is the identity operator and $K^{(\kappa_0)}$ is the boundary integral operator defined by
\[
    \paren{K^{(\kappa_0)} \jb}(\xb) := \nv \times p.v. \int_\Gm \nabla_{\xb} \dfrac{e^{-\iota \kappa_0 R}}{4\pi R} \times \jb(\yb) \ds_{\yb},
\]
where $R = \abs{\xb - \yb}$ and $p.v.$ stands for the Cauchy principal value.

Let the boundary $\Gm$ be discretized into a triangular mesh $\Gm_h$ consisting of $N_f$ triangles, $N_e$ edges, and $N_v$ vertices. We equip on $\Gm_h$ two boundary element spaces $\mathrm{RT}(\Gm_h)$ and $\mathrm{BC}(\Gm_h)$, spanned by the Rao--Wilton--Glisson (RWG) basis functions $\fb_n(\xb)$ \cite{RWG1982} and the BC basis functions $\gb_n(\xb)$ \cite{BC2007}, respectively, with $n = 1, 2, \ldots, N_e$. 

The unknown electric current density $\jb$ is approximated by its expansion in $\mathrm{RT}(\Gm_h)$
\begin{equation}
    \label{eq:expansion}
    \jb(\xb) \approx \sum\limits_{n = 1}^{N_e} \jo_n \fb_n(\xb).
\end{equation}
The coefficient vector is denoted by $\jo := (\jo_1, \jo_2, \ldots, \jo_{N_e})^{\transpose}$. Substituting the expansion \eqref{eq:expansion} into \eqref{eq:MFIE} and testing \eqref{eq:MFIE} with the BC basis functions yields the mixed discretized system of the MFIE\footnote{Plain upright symbols (e.g., $\Go, \Ko, \Po, \Do$) denote matrices associated with the MFIE. Boldface symbols (e.g., $\GGs, \MMs, \PPo, \DDs$) denote matrices associated with the M\"{u}ller formulation, which has a $2\times 2$ operator-block structure. Blackboard-bold symbols (e.g., $\G, \Pbb, \M$) are reserved for matrices in the MT-M\"{u}ller formulation, which has an $N \times N$ subdomain-block structure, with each subdomain block itself having a $2 \times 2$ operator structure.}
\begin{equation}
    \label{eq:dis_MFIE}
    \paren{\dfrac{1}{2} \Go - \Ko^{(\kappa_0)}} \jo = \ho^{inc},
\end{equation}
where the matrices $\Go$ and $\Ko^{(\kappa_0)}$ are defined by
\[
    \left[\Go\right]_{mn} = \inprod{\gb_m, \fb_n}_{\times, \Gm}, \q \left[\Ko^{(\kappa_0)}\right]_{mn} = \inprod{\gb_m, K^{(\kappa_0)} \fb_n}_{\times, \Gm}, 
\]
and the right-hand side vector $\ho^{inc}$ is given by
\[
    \left[\ho^{inc}\right]_m = \inprod{\gb_m, \nv \times \hb^{inc}}_{\times, \Gm}.
\]
Here and throughout the paper, the pairing of two traces supported on $\Gm$ is the skew-symmetric pairing
\[
    \inprod{\gb, \fb}_{\times, \Gm} = \int_{\Gm} (\nv \times \gb) \cdot \fb \ds.
\]

For comparison, we also introduce the classical non-conforming discretization of the MFIE
\begin{equation}
    \label{eq:classical_MFIE}
    \paren{\dfrac{1}{2} \wt{\Go} - \wt{\Ko}^{(\kappa_0)}} \jo = \wt{\ho}^{inc},
\end{equation}
where the tilded matrices and right-hand side vector are obtained by testing \eqref{eq:MFIE} with the rotated RWG basis functions $\nv \times \fb_m$, rather than the BC basis functions $\gb_m$.

The classical discretization \eqref{eq:classical_MFIE} and the mixed discretization \eqref{eq:dis_MFIE} inherit the conditioning behavior of the MFIE. If the wavenumber $\kappa_0$ is sufficiently far from spurious interior resonances and the scatterer $\Om$ is simply connected, both systems are well-conditioned under mesh refinement and in the low-frequency regime, and therefore typically exhibit fast convergence of iterative solvers. However, this favorable behavior is lost in two situations. First, as $\kappa_0$ approaches a resonance, the MFIE operator becomes nearly singular, leading to resonant instability and increasingly ill-conditioned linear systems. Second, on multiply connected surfaces, the static exterior MFIE operator, i.e., the MFIE operator at zero wavenumber, possesses a non-trivial nullspace associated with toroidal loops \cite{CAO+2009a}. Consequently, the mixed discretized system, which correctly preserves the low-frequency behavior of the continuous MFIE, exhibits deteriorating conditioning as the mesh density increases or as the frequency tends to zero. These two sources of deterioration are induced by resonant nullspaces and topology-dependent harmonic nullspaces, and cannot be eliminated by purely algebraic preconditioning.

Despite their similar conditioning properties, the two discretizations differ substantially in accuracy. It is well known that the classical discretization \eqref{eq:classical_MFIE} is far less accurate than the EFIE, particularly in far-field computations \cite{ZCC+2003}. The mixed discretization \eqref{eq:dis_MFIE} significantly improves this accuracy, in some cases bringing it to a level comparable to that of the EFIE \cite{CAD+2011,YJN2011}. Of particular relevance to the present work, \cite{BCA+2014} shows that the mixed MFIE reproduces the correct low-frequency asymptotic scalings of the continuous solution (up to quadrature errors), whereas the classical discretization yields incorrect scalings. Although this correct scaling behavior is not sufficient to guarantee accurate far-field evaluation at low frequencies in finite-precision arithmetic, it provides the necessary foundation for the stabilization schemes developed in \cite{MBC+2020} and in Section~\ref{subsec:MFIE_stabilization}.

\subsection{Low-Frequency Stabilization}
\label{subsec:MFIE_stabilization}

In this section, we propose a stabilization scheme for the mixed MFIE \eqref{eq:dis_MFIE} that suppresses low-frequency cancellation errors while retaining its conditioning properties. The main idea is to separate the solenoidal and non-solenoidal components of the surface current by means of quasi-Helmholtz projectors, and then to rescale these components according to their low-frequency asymptotic scalings.

The non-solenoidal subspace of the trial space $\mathrm{RT}(\Gm_h)$ is conventionally identified with $\range(\wt{\Sm})$, where $\wt{\Sm} \in \R^{N_e \times N_f}$ is the face-edge connectivity matrix defined by \cite{ACB+2013}
\[
    \left[\wt{\Sm}\right]_{mn} :=
    \begin{cases}
        1 \q & \text{if the triangle } n \text{ is identical with } t^+_m, \\
        -1 & \text{if the triangle } n \text{ is identical with } t^-_m, \\
        0 & \text{otherwise}.
    \end{cases}
\]
Here, the RWG basis function $\fb_m(\xb)$ associated with the edge $e_m$ represents a current flowing from the triangle $t_m^+$ to the triangle $t_m^-$. The non-solenoidal subspace has dimension $N_f - 1$. The solenoidal subspace is then defined as the orthogonal complement of the non-solenoidal subspace. Their corresponding projectors on $\mathrm{RT}(\Gm_h)$ are 
\begin{equation}
    \label{eq:qHP}
    \Po^\Sm := \wt{\Sm} \paren{\wt{\Sm}^{\transpose} \wt{\Sm}}^+ \wt{\Sm}^\transpose, \qqq \Po^{\Lambda\textrm{H}} := \Io - \Po^\Sm,
\end{equation}
where $^+$ denotes the Moore--Penrose pseudo-inverse and $\Io$ is the identity matrix \cite{ACB+2013}. Roughly speaking, $\Po^\Sm$ extracts the non-solenoidal component of a coefficient vector in $\mathrm{RT}(\Gm_h)$, whereas $\Po^{\Lambda\textrm{H}}$ extracts its solenoidal component. These projectors are stable with respect to mesh refinement and are mutually orthogonal, i.e., $\Po^\Sm \Po^{\Lambda\textrm{H}} = \Po^{\Lambda\textrm{H}} \Po^\Sm = \zrb$.

Since the mixed discretized MFIE \eqref{eq:dis_MFIE} preserves the low-frequency asymptotic scaling of the continuous solution components at the discrete level \cite{BCA+2011,BCA+2014}, the plane-wave scalings in Table~\ref{tab:scaling} can be expressed for the solution of \eqref{eq:dis_MFIE} in terms of the quasi-Helmholtz projectors as
\begin{equation}
    \label{eq:j_scaling}
    \abs{\Po^\Sm \jo} = \OO(\omega), \qqqqq \abs{\Po^{\Lambda\textrm{H}} \jo} = \OO(1).
\end{equation}
Thus, as $\omega \to 0$, the non-solenoidal component is asymptotically smaller than the solenoidal component. When both components are stored together in finite-precision arithmetic, this imbalance may cause a loss of significant digits in the non-solenoidal component, which in turn can produce inaccurate far-field evaluations.

To balance the two components, we introduce the auxiliary unknown
\[
    \io = \Po^{-1} \jo,
\]
with the right rescaling operator 
\begin{equation}
    \label{eq:right_rescaling}
    \Po = (\iota \omega c_0^{-1} d) \Po^\Sm + \Po^{\Lambda\textrm{H}} .
\end{equation}
Here, $d$ is a characteristic length of the scatterer, introduced to render the factor $\iota \omega c_0^{-1} d$ dimensionless. The imaginary unit $\iota$ is included to ensure the proper recovery of both the real and imaginary parts of the rescaled components, as clarified in \cite{ACB+2013}. It is also convenient for extensions to the time domain, where $\iota\omega$ corresponds to a time derivative.

With this definition, the components of the auxiliary unknown satisfy
\[
    \abs{\Po^\Sm \io} = \OO(1), \qqq \abs{\Po^{\Lambda\textrm{H}} \io} = \OO(1).
\]
The rescaled unknown $\io$ therefore contains components of comparable magnitude. These components can then be stored together in a single floating-point variable without loss of significant digits. 

Substituting $\jo = \Po \io$ into the mixed MFIE \eqref{eq:dis_MFIE} gives
\begin{equation}
    \label{eq:transformed_MFIE}
    \paren{\dfrac{1}{2} \Go - \Ko^{(\kappa_0)}} \Po \io = \ho^{inc}.
\end{equation}
Although the right rescaling provides a cancellation-free representation of the unknown, it does not by itself yield a stable numerical formulation. First, the quasi-Helmholtz components of the right-hand side $\ho^{inc}$ remain asymptotically unbalanced. As a result, the smaller right-hand side component may still suffer from loss of significant digits, which can compromise the accuracy of the corresponding component of $\io$. Second, the right rescaling alters the conditioning of the original mixed MFIE. In particular, the condition number of $\Po$ grows as $\omega^{-1}$ as $\omega\to0$, thereby introducing an artificial source of low-frequency ill-conditioning. A complete stabilization must therefore retain the balanced representation of the unknown while simultaneously balancing the excitation and compensating for the conditioning deterioration caused by the right rescaling. This naturally motivates the application of a left rescaling with asymptotic scalings opposite to those of $\Po$. However, such a rescaling cannot be chosen solely from the scaling of the unknown and excitation. Its construction must also account for the low-frequency behavior of the operator blocks associated with the quasi-Helmholtz components of the trial and test spaces.

To analyze the low-frequency scalings of the quasi-Helmholtz blocks, we decompose the trial and test spaces into their three discrete Helmholtz subspaces. For the trial space $\mathrm{RT}(\Gm_h)$, the local-loop subspace is $\range(\wt{\Lambda})$, where $\wt{\Lambda} \in \R^{N_e \times N_v}$ is the vertex-edge connectivity matrix given by
\[
    \left[\wt{\Lambda}\right]_{mn} :=
    \begin{cases}
        1 \q & \text{if the vertex } n \text{ is identical with } v^+_m, \\
        -1 & \text{if the vertex } n \text{ is identical with } v^-_m, \\
        0 & \text{otherwise}.
    \end{cases}
\]
Here, the edge $e_m$ is oriented from $v^-_m$ to $v^+_m$. The local-loop subspace has dimension $N_v - 1$ \cite{ACB+2013}. The quasi-harmonic subspace is then defined as the orthogonal complement of the non-solenoidal and local-loop subspaces. This subspace has dimension $2g$, where $g$ is the genus of $\Gm$ satisfying 
\begin{equation}
    \label{eq:Euler}
    2g = N_e - N_f - N_v + 2.
\end{equation}
For the test space $\mathrm{BC}(\Gamma_h)$, the roles of the two connectivity matrices are interchanged: $\range(\wt{\Sm})$ represents the local-loop subspace, whereas $\range(\wt{\Lambda})$ represents the non-solenoidal subspace. The quasi-harmonic subspace remains unchanged.

Next, let $\Sm \in \R^{N_e \times (N_f - 1)}$ and $\Lambda \in \R^{N_e \times (N_v - 1)}$ be matrices whose columns form bases for the ranges of $\wt{\Sm}$ and $\wt{\Lambda}$, respectively. Such matrices may be obtained, for example, by removing one linearly dependent column from each connectivity matrix. Let further $\mathrm{H} \in \R^{N_e \times 2g}$ be a matrix whose columns form a basis for the quasi-harmonic subspace. The transformations from the discrete Helmholtz coordinates to the RWG and BC coefficient spaces are then defined, respectively, by
\[
    \Do := \left[\Sm \,\,\, \textrm{H} \,\, \Lambda\right], \qqqq \wt{\Do} := \left[\Lambda \,\,\, \textrm{H} \,\, \Sm\right].
\]
Euler's relation \eqref{eq:Euler} implies that $\Do, \wt{\Do} \in \R^{N_e \times N_e}$ are non-singular. We emphasize that these bases and transformations are introduced solely for the analysis of the low-frequency block scalings. The numerical implementation relies only on the quasi-Helmholtz projectors in \eqref{eq:qHP}.

In the bases of the discrete Helmholtz subspaces, the low-frequency block scalings of the mixed MFIE are known to be \cite{BCA+2011,LCA+2024}
\begin{equation}
    \label{eq:scaling_MFIE}
    \wt{\Do}^\transpose \paren{\dfrac{1}{2} \Go - \Ko^{(\kappa_0)}} \Do = 
    \OO
    \begin{pmatrix}
        1 & \omega^2 & \omega^2 \\
        1 & 1 & \omega^2 \\
        1 & 1 & 1
    \end{pmatrix}.
\end{equation}
A further decomposition of the quasi-harmonic component into poloidal and toroidal loops provides a characterization of the topology-induced nullspaces of the static MFIE operators. However, this finer decomposition is not needed for the present conditioning analysis. We therefore do not pursue it here and refer the reader to \cite{BCA+2011,ADC+2021} for details.

Since the non-solenoidal part of $\Po$ is rescaled by the factor $\iota \omega c_0^{-1} d$ of order $\OO(\omega)$, the transformed system \eqref{eq:transformed_MFIE} has the block scalings
\begin{equation}
    \label{eq:scaling_transformed_MFIE}
    \wt{\Do}^{\transpose} \paren{\dfrac{1}{2} \Go - \Ko^{(\kappa_0)}} \Po \Do =     
    \OO
    \begin{pmatrix}
        \omega & \omega^2 & \omega^2 \\
        \omega & 1 & \omega^2 \\
        \omega & 1 & 1
    \end{pmatrix}.
\end{equation}
This scaling confirms that the right rescaling operator introduces an additional source of low-frequency breakdown. In the dominant lower-triangular part of \eqref{eq:scaling_transformed_MFIE}, all blocks in the first column scale as $\OO(\omega)$, whereas the remaining blocks are of order $\OO(1)$. As a result, the system matrix in \eqref{eq:transformed_MFIE} becomes ill conditioned as $\omega \to 0$, even for simply connected surfaces.

One might attempt to compensate for this ill-conditioning by applying a left rescaling based on the quasi-Helmholtz projectors associated with the $\mathrm{BC}(\Gm_h)$ space. Although this approach restores the desired scaling on simply connected surfaces, it remains problematic on multiply connected geometries. In $\mathrm{BC}(\Gm_h)$, the local-loop and non-solenoidal components exchange roles. Consequently, rescaling the local-loop component in $\mathrm{BC}(\Gm_h)$ also introduces an undesired factor of order $\omega^{-1}$ on the global-loop component. The condition number may therefore still grow as $\OO(\omega^{-1})$ as $\omega \to 0$. This is analogous to the difficulty identified in \cite{BMC+2017} for the PMCHWT formulation.

To avoid this problem, we adapt the coefficient-space rescaling strategy of \cite{BMC+2017}. We first multiply \eqref{eq:transformed_MFIE} by the inverse Gram matrix $\Go^{-1}$, obtaining
\begin{equation}
    \label{eq:transformed_MFIE_1}
    \paren{\dfrac{1}{2} \Io - \Go^{-1} \Ko^{(\kappa_0)}} \Po \io = \Go^{-1} \ho^{inc}.
\end{equation}
The Gram matrix $\Go$ is invertible on closed surfaces. The key effect of this multiplication is that the resulting operator maps the $\mathrm{RT}(\Gm_h)$ coefficient space into itself. Thus, both the right and left rescalings can be applied in the same coefficient space.

The scaling of \eqref{eq:transformed_MFIE_1} follows from the algebraic structure
\begin{equation}
    \label{eq:Go_inv}
    \Do^{-1} \Go^{-1} \wt{\Do}^{-\transpose} = 
    \begin{pmatrix}
        \square & 0 & 0 \\
        \square & \square & 0 \\
        \square & \square & \square
    \end{pmatrix},
\end{equation}
where $\square$ denotes a non-zero block \cite{BCA2015b}. Hence,
\begin{equation}
    \label{eq:scaling_transformed_MFIE_1}
    \Do^{-1} \paren{\dfrac{1}{2} \Io - \Go^{-1} \Ko^{(\kappa_0)}} \Po \Do =     
    \OO
    \begin{pmatrix}
        \omega & \omega^2 & \omega^2 \\
        \omega & 1 & \omega^2 \\
        \omega & 1 & 1
    \end{pmatrix}.
\end{equation}
Although the block scalings in \eqref{eq:scaling_transformed_MFIE_1} appear similar to those in \eqref{eq:scaling_transformed_MFIE}, their interpretation is different. In \eqref{eq:scaling_transformed_MFIE}, the rows represent dual components associated with the test space $\mathrm{BC}(\Gamma_h)$. In contrast, the rows in \eqref{eq:scaling_transformed_MFIE_1} represent components in the $\mathrm{RT}(\Gamma_h)$ coefficient space. This change of range is reflected by the appearance of $\wt{\Do}^{\transpose}$ in \eqref{eq:scaling_transformed_MFIE} and $\Do^{-1}$ in \eqref{eq:scaling_transformed_MFIE_1}.

We can now apply the left rescaling in the $\mathrm{RT}(\Gamma_h)$ coefficient space. Left multiplying \eqref{eq:transformed_MFIE_1} by
\[
     \Po^{-1} = (\iota \omega c_0^{-1} d)^{-1} \Po^\Sm + \Po^{\Lambda\textrm{H}}
\]
gives the stabilized MFIE system
\begin{equation}
    \label{eq:stabilized_MFIE}
    \paren{\dfrac{1}{2} \Io - \Po^{-1} \Go^{-1} \Ko^{(\kappa_0)} \Po} \io = \Po^{-1} \Go^{-1} \ho^{inc}.
\end{equation}
With respect to the discrete Helmholtz decomposition, the system matrix has the block scalings
\[
    \Do^{-1} \paren{\dfrac{1}{2} \Io - \Po^{-1} \Go^{-1} \Ko^{(\kappa_0)} \Po} \Do =     
    \OO
    \begin{pmatrix}
        1 & \omega & \omega \\
        \omega & 1 & \omega^2 \\
        \omega & 1 & 1
    \end{pmatrix},
\]
whereas the right-hand side satisfies \cite{BCA+2011}
\[
    \abs{\Po^\Sm \Po^{-1} \Go^{-1} \ho^{inc}} = \OO(1), \q\, \abs{\Po^{\Lambda\textrm{H}} \Po^{-1} \Go^{-1} \ho^{inc}} = \OO(1).
\]
The left rescaling therefore compensates for the low-frequency breakdown introduced by the right rescaling and simultaneously balances the right-hand side. Consequently, \eqref{eq:stabilized_MFIE} is free from low-frequency cancellation errors while retaining the conditioning behavior of the mixed MFIE. In particular, for simply connected geometries and at frequencies sufficiently far from interior resonances, the system remains well conditioned under mesh refinement and in the low-frequency regime. Nevertheless, it inherits the intrinsic resonant instability of the MFIE and, for multiply connected geometries, the topology-induced breakdown at low frequencies. However, the numerical results in Section~\ref{sec:results} show that the topology-induced breakdown does not prevent rapid solution by means of Krylov iterative solvers and does not compromise the accuracy of the computed far field. 

A comprehensive treatment of both topology-induced breakdown and resonant instability was proposed in \cite{MBC+2020}. In that approach, the exterior MFIE is combined with its interior counterpart to obtain a symmetrized formulation. Together with a quasi-Helmholtz decomposition, the symmetrized operator structure gives access to all the blocks governed by static cancellation identities and those associated with the nullspaces of the static MFIE operators on multiply connected geometries. The corresponding static contributions can therefore be set exactly to zero, preventing quadrature errors from destroying these cancellations or perturbing the nullspace. When combined with a preconditioned EFIE, the formulation yields a stabilized CFIE that is free from low-frequency and dense-mesh breakdown, as well as topology-induced and resonant instabilities, while preserving the exact cancellations and nullspace structure of the static MFIE operators.

The formulation proposed here is structurally simpler than the symmetrized approach of \cite{MBC+2020}, as it avoids the additional matrix-vector products required by the symmetrized operator. More importantly, it extends naturally to the M\"uller equation, whose original mixed formulation is not inherently affected by topology-induced breakdown or resonant instability. By contrast, extending the symmetrized MFIE approach would require an analogue of the interior MFIE operator, and no such operator is readily available for the M\"uller formulation.

Finally, after solving the stabilized system \eqref{eq:stabilized_MFIE} for $\io$, the coefficient vector $\jo$ of the physical electric current density can be obtained by applying the rescaling operator $\Po$ to $\io$. To avoid reintroducing cancellation errors, the corresponding far-field pattern should be evaluated directly from $\io$ as follows:
\begin{equation}
    \label{eq:dis_farfield}
    \begin{aligned}
        \Eb^{\text{far}}(\wh{\xb}) & = (\iota \omega c_0^{-1} d) \sum\limits_{n = 1}^{N_e} \left[\Po^{\Sm} \io \right]_n \int_\Gm e^{-\iota \kappa_0 \, \wh{\xb} \cdot \yb} \, \fb_n(\yb) \ds_{\yb} \\
        & + \sum\limits_{n = 1}^{N_e} \left[\Po^{\Lambda\textrm{H}} \io \right]_n \int_\Gm \paren{e^{-\iota \kappa_0 \, \wh{\xb} \cdot \yb} - 1} \fb_n(\yb)\ds_{\yb}.
    \end{aligned}
\end{equation}
Here, $\left[\Po^{\Sm} \io \right]_n$ and $\left[\Po^{\Lambda\mathrm{H}}\io\right]_n$ denote the $n$-th entries of the projected coefficient vectors $\Po^{\Sm}\io$ and $\Po^{\Lambda\mathrm{H}}\io$, respectively.

\section{M\"uller Integral Equation}
\label{sec:Mueller}

In this section, we extend the low-frequency regularization developed for the MFIE to the M\"{u}ller equation for homogeneous dielectric objects. The extension is natural because the two formulations share the same essential structure: both are second-kind boundary integral equations and both admit accurate mixed discretizations based on dual testing. 

\subsection{Formulation and Discretization}
\label{subsec:Mueller_formulation}

We adopt the setting described in Section~\ref{subsec:MFIE_formulation} and consider scattering by a homogeneous dielectric object $\Om$ with permittivity $\epsilon_1$ and permeability $\mu_1$. The magnetic and electric surface current densities $\mb$ and $\jb$ induced on $\Gm$ satisfy the M\"uller integral equation \eqref{eq:Mueller} \cite{Muller1969}. In this formulation, the single-layer operator $T^{(\kappa_0)}$ associated with the background medium is defined by
\begin{align*}
    \paren{T^{(\kappa_0)} \jb} (\xb) & := -\iota \kappa_0 \nv \times \int_{\Gm} \dfrac{e^{-\iota \kappa_0 R}}{4\pi R} \jb(\yb) \ds_{\yb} \\
    & \q\,\,\, + \dfrac{1}{\iota \kappa_0} \nv \times \nabla \int_{\Gm} \dfrac{e^{-\iota \kappa_0 R}}{4\pi R} \nabla_\Gm \cdot \jb(\yb) \ds_{\yb}.
\end{align*}
The operators $T^{(\kappa_1)}$ and $K^{(\kappa_1)}$ associated with the interior medium are defined analogously to $T^{(\kappa_0)}$ and $K^{(\kappa_0)}$, with wavenumber $\kappa_1 = \omega/c_1$ and wave speed $c_1 = 1/\sqrt{\mu_1 \epsilon_1}$.

\begin{figure*}[!t]
    \begin{align}
        \label{eq:Mueller}
        \begin{pmatrix}
            \dfrac{\epsilon_1 + \epsilon_0}{2} I + \epsilon_1 K^{(\kappa_1)} - \epsilon_0 K^{(\kappa_0)} & - c_1^{-1} T^{(\kappa_1)} + c_0^{-1} T^{(\kappa_0)} \\
            c_1^{-1} T^{(\kappa_1)} - c_0^{-1} T^{(\kappa_0)} & \dfrac{\mu_1 + \mu_0}{2} I + \mu_1 K^{(\kappa_1)} - \mu_0 K^{(\kappa_0)}
        \end{pmatrix}
        \begin{pmatrix}
            \mb \\[0.2cm]
            \jb
        \end{pmatrix}
         & = 
         \begin{pmatrix}
             \epsilon_0 \, \eb^{inc} \times \nv \\[0.2cm]
             \mu_0 \, \nv \times \hb^{inc}
         \end{pmatrix}.
         \\[0.2cm]
        \label{eq:dis_Mueller}
        \begin{pmatrix}
            \dfrac{\epsilon_1 + \epsilon_0}{2} \Go + \epsilon_1 \Ko^{(\kappa_1)} - \epsilon_0 \Ko^{(\kappa_0)} & - c_1^{-1} \To^{(\kappa_1)} + c_0^{-1} \To^{(\kappa_0)} \\
            c_1^{-1} \To^{(\kappa_1)} - \, c_0^{-1} \To^{(\kappa_0)} & \dfrac{\mu_1 + \mu_0}{2} \Go + \mu_1 \Ko^{(\kappa_1)} - \mu_0 \Ko^{(\kappa_0)}
        \end{pmatrix}
        \begin{pmatrix}
            \mo \\[0.2cm]
            \jo
        \end{pmatrix}
         & = 
        \begin{pmatrix}
            \epsilon_0 \eo^{inc} \\[0.2cm]
            \mu_0 \ho^{inc}
        \end{pmatrix}.
    \end{align}
    \hrulefill
\end{figure*}

In the M\"uller equation, the exterior and interior representation formulas are combined so that the hyper-singularities of the single-layer operators appearing in the off-diagonal blocks cancel. This hyper-singularity cancellation is the defining feature of the M\"uller formulation, which leads to well-conditioned discretized matrices in both the dense-mesh and low-frequency regimes \cite{BCA+2014,LC2026}.

The unknown current pair $(\mb, \jb)^\transpose$ is approximated in the product trial space $\bm{\mathrm{RT}}(\Gm_h) := \mathrm{RT}(\Gm_h) \times \mathrm{RT}(\Gm_h)$, and its coefficient vector is denoted by $\uv = (\mo, \jo)^\transpose$. The M\"uller equation \eqref{eq:Mueller} is tested with the BC basis functions in the product test space $\bm{\mathrm{BC}}(\Gm_h) := \mathrm{BC}(\Gm_h) \times \mathrm{BC}(\Gm_h)$, yielding the mixed discretization
\begin{equation}
    \label{eq:dis_Mueller_1}
    \MMs \uv = \uv^{inc}.
\end{equation}
Its explicit form is given in \eqref{eq:dis_Mueller}. In particular, the matrix $\To^{(\kappa_0)}$ and the right-hand side vector $\eo^{inc}$ are defined by
\begin{align*}
    \left[\To^{(\kappa_0)}\right]_{mn} & = \inprod{\gb_m, T^{(\kappa_0)} \fb_n}_{\times, \Gm}, \\
    \left[\eo^{inc}\right]_{m} & = \inprod{\gb_m, \eb^{inc} \times \nv}_{\times, \Gm}.
\end{align*}
The matrices $\Go$ and $\Ko^{(\kappa_0)}$, together with the vector $\ho^{inc}$, have already been defined in \eqref{eq:dis_MFIE}. The matrices $\To^{(\kappa_1)}$ and $\Ko^{(\kappa_1)}$ are defined analogously to $\To^{(\kappa_0)}$ and $\Ko^{(\kappa_0)}$.

For comparison, we also introduce the classical non-conforming discretization of the M\"uller equation
\begin{equation}
    \label{eq:classical_Mueller}
     \wt{\MMs} \, \uv = \wt{\uv}^{inc},
\end{equation}
where the corresponding tilded matrices and vectors are obtained by testing \eqref{eq:Mueller} with the rotated RWG basis functions instead of the BC basis functions.

Both the classical and mixed discretizations inherit the well-conditioning of the M\"uller equation. In particular, unlike the MFIE, the M\"{u}ller equation remains well conditioned in the dense-mesh and low-frequency regimes on both simply and multiply connected geometries, and it does not suffer from spurious resonances. As in the MFIE case, the mixed discretization generally yields more accurate solutions than the classical discretization \cite{YJN2011,YJN2013,YKJ2016}. More importantly, it also preserves the low-frequency asymptotic scalings of the Helmholtz components of the continuous solution \cite{BCA+2014}. This property is essential for the stabilization developed in Section~\ref{subsec:Mueller_stabilization}, but is guaranteed neither by the classical discretization \eqref{eq:classical_Mueller} nor by the other non-conforming or discontinuous Galerkin schemes considered in \cite{YT2005,UTR2011,CHS2017}.

\subsection{Low-Frequency Stabilization}
\label{subsec:Mueller_stabilization}

Since the M\"{u}ller formulation is a $2 \times 2$ block system, we first introduce the corresponding block-diagonal extensions of the matrices and projectors introduced in Section~\ref{subsec:MFIE_stabilization}
\begin{align*}
    & \GGs = \diag(\Go, \Go), && \PPo = \diag(\Po, \Po), \\
    & \PPo^{\Sm} = \diag(\Po^\Sm, \Po^\Sm), && \PPo^{\Lambda\textrm{H}} = \diag(\Po^{\Lambda\textrm{H}}, \Po^{\Lambda\textrm{H}}), \\
    & \DDs = \diag(\Do, \Do), && \wt{\DDs} = \diag(\wt{\Do}, \wt{\Do}).
\end{align*}

The wave speed appearing in the rescaling factor of $\Po$ may be chosen as either $c_0$ or $c_1$. This choice is not critical, since both are independent of $\omega$ and therefore provide the same cancellation-error mitigation. Moreover, numerical results demonstrate that the conditioning of the resulting formulations is essentially unaffected by this choice, even in high-contrast configurations.

As in the MFIE case, the idea is to rescale the quasi-Helmholtz components of the surface currents so that their magnitudes remain balanced as $\omega \to 0$. For dielectric scattering, the components of both the electric and magnetic surface currents exhibit the low-frequency scaling exponents
\cite{GAM+2017,BMC+2017,LMA+2026}
\[
    \alpha = \beta = 1, \qqqqq \gm = 0.
\]
Thus, unlike in PEC scattering, the global- and local-loop components scale differently. The quasi-Helmholtz projectors employed in Section~\ref{subsec:MFIE_stabilization} separate only the non-solenoidal and solenoidal subspaces and therefore cannot rescale the global- and local-loop components independently. As a result, the global-loop component, after rescaling, remains susceptible to finite-precision cancellation. Nevertheless, its contribution to the far field is asymptotically negligible compared with that of the local-loop component, as follows from \eqref{eq:farfield}. Recovering the global-loop component would therefore provide no substantial improvement in far-field accuracy. Consequently, it is sufficient to separate and rescale the non-solenoidal and combined solenoidal components. No further decomposition of the solenoidal subspace into global and local loops is required.

Using the solenoidal and non-solenoidal quasi-Helmholtz projectors, the low-frequency asymptotic scaling of the solution $\uv$ of \eqref{eq:dis_Mueller_1} can be expressed as
\[
    \abs{\PPo^\Sm \uv} = \OO(\omega), \qqqq \abs{\PPo^{\Lambda\textrm{H}} \uv} = \OO(1),
\]
which is the same component-wise scaling as in \eqref{eq:j_scaling}. The physical unknown is therefore rescaled by introducing the auxiliary unknown
\[
    \vvs = \PPo^{-1} \uv,
\]
or equivalently $\uv = \PPo \vvs$. By construction, the components of $\vvs$ satisfy
\[
    \abs{\PPo^\Sm \vvs} = \OO(1), \qqqq \abs{\PPo^{\Lambda\textrm{H}} \vvs} = \OO(1).
\]
Thus, expressing the mixed M\"uller system in terms of $\vvs$ eliminates the low-frequency imbalance between the quasi-Helmholtz components of the physical surface currents. This balanced representation is the key ingredient for evaluating the far field accurately in the low-frequency regime.

Next, we adopt the coefficient-space rescaling strategy developed for the mixed MFIE in Section~\ref{subsec:MFIE_stabilization}. The right rescaling by $\PPo$ yields a cancellation-free representation of the unknowns but introduces an artificial low-frequency breakdown in the system matrix. Moreover, the right-hand side remains unbalanced. To compensate for these effects, we left-multiply the system by the inverse Gram matrix $\GGs^{-1}$ and subsequently apply the inverse rescaling $\PPo^{-1}$. The resulting stabilized M\"uller formulation reads
\begin{equation}
   \label{eq:stabilized_Mueller}
   \PPo^{-1} \GGs^{-1} \MMs \PPo \vvs = \PPo^{-1} \GGs^{-1} \uv^{inc}.
\end{equation}

To justify this construction, we first record the low-frequency block scalings of the mixed M\"{u}ller matrix in \eqref{eq:dis_Mueller_1}. With respect to the quasi-Helmholtz decompositions of the trial and test spaces, one obtains
\[
    \wt{\DDs}^\transpose \MMs \, \DDs = 
    \OO\paren{
        \begin{array}{c c c | c c c}
            1 & \omega^2 & \omega^2 & \omega & \omega & \omega\\
            1 & 1 & \omega^2 & \omega & \omega & \omega\\
            1 & 1 & 1 & \omega & \omega & \omega \\ \hline        
            \omega & \omega & \omega & 1 & \omega^2 & \omega^2 \\
            \omega & \omega & \omega & 1 & 1 & \omega^2 \\
            \omega & \omega & \omega & 1 & 1 & 1
        \end{array}
    }.
\]
The diagonal operator blocks have the same structure as the mixed MFIE and therefore inherit the corresponding block scalings. The off-diagonal operator blocks contain the single-layer operators $T^{(\kappa_0)}$ and $T^{(\kappa_1)}$. Their weakly singular (vector-potential) contributions scale as $\OO(\omega)$, whereas the hyper-singular (scalar-potential) contributions scale individually as $\OO(\omega^{-1})$. In the specific linear combination in the M\"{u}ller equation, however, the leading hyper-singular terms cancel, leaving contributions of order $\OO(\omega)$ \cite{BCA+2014,LC2026}. Consequently, all Helmholtz sub-blocks inside the off-diagonal operator blocks have scaling $\OO(\omega)$.

A direct scaling calculation for the matrix in the stabilized formulation \eqref{eq:stabilized_Mueller} gives
\[
    \DDs^{-1} \PPo^{-1} \GGs^{-1} \MMs \PPo \DDs = 
    \OO\paren{
        \begin{array}{c c c | c c c}
            1 & \omega & \omega & \omega & 1 & 1 \\
            \omega & 1 & \omega^2 & \omega^2 & \omega & \omega\\
            \omega & 1 & 1 & \omega^2 & \omega & \omega \\ \hline        
            \omega & 1 & 1 & 1 & \omega & \omega \\
            \omega^2 & \omega & \omega & \omega & 1 & \omega^2 \\
            \omega^2 & \omega & \omega & \omega & 1 & 1
        \end{array}
    }.
\]
This shows that the left rescaling compensates for the artificial low-frequency ill-conditioning introduced by the right rescaling. The resulting right-hand side can also be verified to have balanced component-wise scaling. Consequently, \eqref{eq:stabilized_Mueller} eliminates cancellation errors arising from the loss of significant digits in both the unknown and right-hand side components, while retaining the favorable conditioning of the original mixed M\"uller system. In particular, the formulation remains well conditioned under mesh refinement and in the low-frequency regime for both simply and multiply connected geometries, and it does not suffer from resonant instability.

After solving \eqref{eq:stabilized_Mueller} for $\vvs$, the coefficient vectors of the physical magnetic and electric surface currents are recovered by $\uv = \PPo \vvs$.  For far-field evaluations, however, the cancellation-free representation $\vvs$ is used directly, rather than first forming the physical unknown $\uv$, thereby ensuring accuracy in the low-frequency regime. The computation is analogous to \eqref{eq:dis_farfield}.

\section{Multitrace M\"uller Equation}
\label{sec:MT_Mueller}


We now extend the proposed low-frequency regularization to the M\"{u}ller equation for scattering by piecewise homogeneous dielectric objects. For this purpose, we employ the recently developed global multitrace (MT) formulation of \cite{LC2026}. This formulation is composed of second-kind operators with M\"{u}ller-type structure and therefore inherits the favorable conditioning properties of the M\"{u}ller equation. It also admits a mixed RWG--BC discretization, which is particularly important in the present context. First, it yields accurate surface-current solutions and derived far field. Second, the dedicated low-frequency scaling analysis of Bogaert et. al. \cite{BCA+2014} applies directly to this global MT mixed discretization, showing that the discrete solution preserves the component-wise low-frequency scaling of its continuous counterpart. Consequently, the regularization strategy developed in Section~\ref{subsec:Mueller_stabilization} can be extended straightforwardly to the MT-M\"uller formulation. 

To the best of our knowledge, the single-trace formulation of \cite{CHS2017} is the only other M\"{u}ller-type formulation for composite penetrable objects. However, its discontinuous Galerkin discretization does not preserve the asymptotic scaling of either the continuous solution components or the Helmholtz sub-blocks of the operator. In addition, a mixed RWG--BC discretization is not available in the single-trace setting for composite geometries containing junctions. This limitation is structural: such a discretization would require a boundary element space that is both dual to the RWG space and a discrete subspace of the single-trace energy space, so that the transmission conditions are enforced. No such space is currently available for junction configurations. The same difficulty arises in the preconditioning of the single-trace PMCHWT formulation for composite objects, which likewise requires an appropriate dual conforming subspace. Although a recently proposed quasi-local multitrace approach provides an alternative remedy \cite{Cools2026}, it is not directly compatible with the regularization framework considered here.

\begin{figure}
    \centering
    \includegraphics[width=\linewidth]{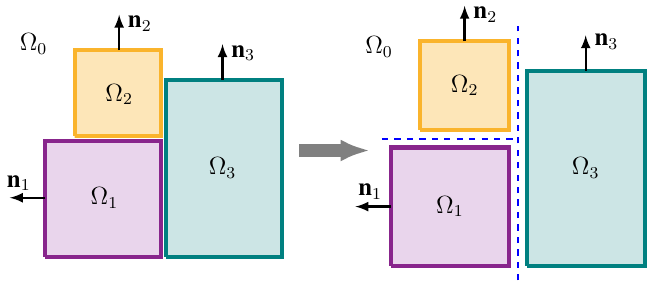}
    \caption{\textit{Left}: A composite object consisting of $N$ components $\Om_k, k = 1, 2, \ldots, N,$ immersed in the background $\Om_0$. \textit{Right}: Illustration of the global multitrace approach, in which a conceptual gap filled by the background medium is inserted between adjacent components.}
    \label{fig:domain}
\end{figure}

We consider the dielectric scattering setting described in Section~\ref{subsec:Mueller_formulation}, where the scatterer $\Om$ is composed of $N$ homogeneous dielectric subdomains
\[
    \ovl{\Om} = \bigcup_{k=1}^N \ovl{\Om_k}.
\]
Each subdomain $\Om_k$ is occupied by a dielectric material with permittivity $\epsilon_k$ and permeability $\mu_k$. Its boundary is denoted by $\Gm_k$, with outward unit normal $\nv_k$ (see Fig.~\ref{fig:domain}, left).

In the global MT framework, the magnetic and electric surface current densities are sought on the boundary of each subdomain $\Om_k$. On exterior interfaces, which belong to the boundary of the background domain $\Om_0$, each current appears only once and the transmission conditions are enforced directly. On interior interfaces shared by adjacent subdomains, however, the currents are duplicated, with one trace assigned to each side of the interface. The transmission conditions on these interior interfaces are then enforced indirectly through the introduction of an infinitesimally thin gap. In this interpretation, all $N$ subdomains are treated as floating in the background medium $\Om_0$ \cite{CHJ2013} (see Fig.~\ref{fig:domain}, right).

We recall the global MT-M\"{u}ller equation introduced in \cite{LC2026}. Let $\ub = \paren{\ub_1, \ub_2, \ldots, \ub_N}^\transpose$, where $\ub_k = \paren{\mb_k, \jb_k}^\transpose$ is the surface currents on the boundary $\Gm_k$ of $\Om_k$. The global current vector $\ub$ satisfies 
\begin{equation}
    \label{eq:MT_Mueller}
    \Mb \ub = \ub^{inc},
\end{equation}
where $\Mb$ has an $N \times N$ subdomain-block structure. Its diagonal and off-diagonal subdomain blocks are defined in \eqref{eq:Mueller_diagonal} and \eqref{eq:Mueller_off_diagonal}, respectively. The subscripts $j$ and $k$ on the boundary integral operators indicate the target and source boundaries, respectively.
The wavenumber and wave speed in $\Om_k$ are $\kappa_k = \omega/c_k$ and $c_k = 1/\sqrt{\mu_k \epsilon_k}$, respectively. The right-hand side is $\ub^{inc} = \paren{\ub^{inc}_1, \ub^{inc}_2, \ldots, \ub^{inc}_N}^\transpose$, where $\ub_k^{inc} = \paren{\epsilon_0 \, \eb^{inc} \times \nv_k, \mu_0 \, \nv_k \times \hb^{inc}}^\transpose$ denotes the trace of the incident fields on $\Gm_k$. 
\begin{figure*}
    \centering
    \begin{align}
    \label{eq:Mueller_diagonal}
        \left[\Mb\right]_{kk} & = 
        \begin{pmatrix}
            \dfrac{\epsilon_0 + \epsilon_k}{2} I + \epsilon_k K_{kk}^{(\kappa_k)} - \epsilon_0 K_{kk}^{(\kappa_0)} & - c_k^{-1} T_{kk}^{(\kappa_k)} + c_0^{-1} T_{kk}^{(\kappa_0)} \\
            c_k^{-1} T_{kk}^{(\kappa_k)} - c_0^{-1} T_{kk}^{(\kappa_0)} & \dfrac{\mu_0 + \mu_k}{2} I + \mu_k K_{kk}^{(\kappa_k)} - \mu_0 K_{kk}^{(\kappa_0)}
        \end{pmatrix}, \q k = 1, 2, \ldots, N, \\[0.2cm]
    \label{eq:Mueller_off_diagonal}
        \left[\Mb\right]_{jk} & = 
        \begin{pmatrix}
            \dfrac{\epsilon_0 - \epsilon_k}{2} I + \epsilon_k K^{(\kappa_k)}_{jk} - \epsilon_0 K^{(\kappa_0)}_{jk} & -c_k^{-1} T^{(\kappa_k)}_{jk} + c_0^{-1} T^{(\kappa_0)}_{jk} \\
            c_k^{-1} T^{(\kappa_k)}_{jk} - c_0^{-1} T^{(\kappa_0)}_{jk} & \dfrac{\mu_0 - \mu_k}{2} I + \mu_k K^{(\kappa_k)}_{jk} - \mu_0 K^{(\kappa_0)}_{jk}
        \end{pmatrix}, \q j, k = 1, 2, \ldots, N, \,\, j \neq k.
    \end{align}   
    \hrulefill
\end{figure*}

Each off-diagonal subdomain block in the $k$-th block column differs from the corresponding diagonal block in two respects. First, the operator $I$ appearing in the off-diagonal subdomain block represents a geometric identity, as explained in \cite{LBG+2022}. Second, the coefficients multiplying the geometric identities are the differences $\epsilon_0 - \epsilon_k$ and $\mu_0 - \mu_k$, rather than the corresponding sums. Both the diagonal and off-diagonal subdomain blocks have a M\"uller-type structure, in which the hyper-singularities of the single-layer operators cancel.

Let the boundary $\Gm_k$ of each subdomain $\Om_k$ be discretized by a triangular mesh $\Gm_{h, k}$. The global mesh is defined as the disjoint union 
\[
    \Gm_h := \bigsqcup_{k = 1}^N \Gm_{h, k}.
\]
In this setting, the local meshes $\Gamma_{h, k}$ and $\Gamma_{h, j}$ are not required to be mutually conformal on the shared interface $\Gamma_j \cap \Gamma_k$. This is one of the main advantages of the MT framework: the subdomain boundaries can be discretized independently, thereby enabling the efficient numerical treatment of geometrically complex and multiscale problems.

On $\Gm_h$, we introduce the product spaces
\[
    \mathbb{RT}(\Gm_h) := \prod_{k=1}^N \mathbf{RT}(\Gm_{h, k}), \qq \mathbb{BC}(\Gm_h) := \prod_{k=1}^N \mathbf{BC}(\Gm_{h, k}), 
\]
where the local spaces $\mathbf{RT}(\Gm_{h, k}) = \mathrm{RT}(\Gm_{h, k}) \times \mathrm{RT}(\Gm_{h, k})$ and $\mathbf{BC}(\Gm_{h, k}) = \mathrm{BC}(\Gm_{h, k}) \times \mathrm{BC}(\Gm_{h, k})$. Here, $\mathrm{RT}(\Gm_{h,k})$ and $\mathrm{BC}(\Gm_{h,k})$ denote the RWG and BC spaces defined on the local mesh $\Gm_{h,k}$, respectively. Note that, even when the local meshes are mutually conformal, the space $\mathbb{BC}(\Gm_h)$ cannot represent single-trace functions satisfying the transmission conditions, unlike $\mathbb{RT}(\Gm_{h})$.

The unknown $\ub$ is approximated in $\mathbb{RT}(\Gm_h)$, and the corresponding coefficient vector is denoted by $\ufrakbm$. Testing the MT-M\"{u}ller equation \eqref{eq:MT_Mueller} with basis functions in $\mathbb{BC}(\Gm_h)$ gives the mixed system
\begin{equation}
    \label{eq:dis_MT_Mueller}
    \M \, \ufrakbm = \ufrakbm^{inc}.
\end{equation}

The mixed MT-M\"{u}ller system \eqref{eq:dis_MT_Mueller} inherits the conditioning properties of the mixed M\"{u}ller equation \eqref{eq:dis_Mueller_1}. In particular, it remains well conditioned in the dense-mesh and low-frequency regimes on both simply and multiply connected geometries, and does not suffer from resonant instability. In addition, the quasi-Helmholtz components of the surface currents on each subdomain exhibit the same low-frequency scalings as those in the single-object formulation.

This analogy makes the extension of the low-frequency stabilization straightforward. We define the subdomain-block diagonal matrices
\begin{align*}
    & \G = \diag(\GGs_1, \GGs_2, \ldots, \GGs_N), \\
    & \Pbb = \diag(\PPo_1, \PPo_2, \ldots, \PPo_N), 
\end{align*}
where $\GGs_k = \diag(\Go_k, \Go_k)$ and $\PPo_k = \diag(\Po_k, \Po_k)$. Here, $\Go_k$ is the Gram matrix coupling $\mathrm{RT}(\Gm_{h,k})$ and $\mathrm{BC}(\Gm_{h,k})$, while $\Po_k$ is the quasi-Helmholtz rescaling operator associated with $\mathrm{RT}(\Gm_{h,k})$. Introducing the auxiliary unknown
\[
    \vfrakbm = \Pbb^{-1} \, \ufrakbm,
\]
and applying the same rescaling procedure as in Section~\ref{subsec:Mueller_stabilization}, we obtain the stabilized MT-M\"{u}ller formulation
\begin{equation}
    \label{eq:stablized_MT_Mueller}
    \Pbb^{-1} \G^{-1} \M \, \Pbb \, \vfrakbm = \Pbb^{-1} \G^{-1} \ufrakbm^{inc}.
\end{equation}

Analogously to the stabilized M\"{u}ller formulation \eqref{eq:stabilized_Mueller}, the system \eqref{eq:stablized_MT_Mueller} is free from low-frequency cancellation errors and preserves the favorable conditioning of the original mixed MT-M\"{u}ller system. It is therefore stable in the dense-mesh and low-frequency regimes, including configurations with non-trivial topology, and it does not suffer from resonant instability.

In composite configurations consisting of multiple subdomains, several choices are available for the rescaling parameters entering $\Pbb$. Each local operator $\PPo_k$ may be defined using parameters associated with the corresponding subdomain $\Om_k$, such as the local wave speed $c_k$ and a characteristic length $d_k$. Alternatively, a uniform rescaling may be adopted using the background wave speed $c_0$ and the diameter $d$ of the entire scatterer $\Om$. For moderate material contrasts and subdomains of comparable size, these alternatives yield similar scaling and conditioning properties. For high-contrast or multiscale configurations, we recommend employing the local parameters $c_k$ and $d_k$. Indeed, the quantity $c_k^{-1} d_k$ appearing in the rescaling operator $\PPo_k$ represents the characteristic wave-travel time across the subdomain $\Om_k$. Consequently, this choice provides a more balanced normalization by reflecting the local material properties and geometric scales.

\section{Implementation Details}
\label{sec:implementation}

This section discusses the implementation of the proposed regularizations for the MFIE and the M\"uller equations, with particular emphasis on the numerical treatment required to preserve the theoretical low-frequency behavior in practical computations. The additional computational cost introduced by the regularization is also analyzed and compared with the asymptotic complexity of the original formulations.

Throughout this section, we consider an incident plane wave of the form
\begin{equation}
    \label{eq:plane_wave}
    \begin{aligned}
        \eb^{inc}(\xb) & =  \pb \, e^{-\iota \kappa_0 \wh{\kb} \cdot \xb}, \\
        \hb^{inc}(\xb) & = \sqrt{\dfrac{\epsilon_0}{\mu_0}} (\wh{\kb} \times \pb) \, e^{-\iota \kappa_0 \wh{\kb} \cdot \xb},
    \end{aligned}
\end{equation}
where $\pb$ denotes the polarization vector and $\wh{\kb}$ is the unit propagation direction, with $\wh{\kb} \cdot \pb = 0$.

\subsection{MFIE}
\label{subsec:MFIE_implementation}

The low-frequency stability of the regularized MFIE \eqref{eq:stabilized_MFIE} relies on a cancellation property of the static MFIE operator. Specifically,
\begin{equation}
    \inprod{\gb, \paren{\dfrac{1}{2} I - K^{(0)}} \fb}_{\times, \Gm} = 0, \label{eq:vanishing_M0}
\end{equation}
when $\fb$ and $\gb$ are both local-loop functions, or when one is a local-loop function and the other is harmonic. Detailed derivations of \eqref{eq:vanishing_M0} are provided in \cite{BCA+2011,ADC+2021}. The role of this identity becomes apparent from the low-frequency expansion of the gradient of the Green function
\[
    \nabla_{\xb} \dfrac{e^{-\iota \kappa_0 R}}{4\pi R} = -\dfrac{\xb - \yb}{4\pi R^3} \paren{1 + \dfrac{\kappa_0^2 R^2}{2} - \dfrac{\iota \kappa_0^3 R^3}{3} + \cdots}.
\]
The absence of a term linear in $\kappa_0$ implies that the leading frequency-dependent term is of order $\kappa_0^2$. Moreover, the static contribution is eliminated in the local-loop/local-loop and local-loop/harmonic pairings by virtue of \eqref{eq:vanishing_M0}. Consequently, the corresponding Helmholtz blocks in \eqref{eq:scaling_MFIE} scale as $\OO(\omega^2)$, whereas the remaining blocks are of order $\OO(1)$. 

This asymptotic block structure is essential to the proposed regularization. The $\OO(\omega^2)$ scaling of the relevant blocks allows the left rescaling to compensate for the low-frequency ill-conditioning introduced by the right rescaling. As a result, the regularized system retains the low-frequency conditioning behavior of the original mixed MFIE.

The preceding argument, however, is based on exact operator identities. At the discrete level, quadrature errors and the non-conformity of the quasi-Helmholtz projectors with respect to the continuous Helmholtz decomposition may prevent the cancellation in \eqref{eq:vanishing_M0} from being reproduced with sufficient accuracy. To identify the block affected by this loss of cancellation, consider the following decomposition of the second term in \eqref{eq:stabilized_MFIE}:
\begin{align*}
    \Po^{-1} \Go^{-1} \Ko^{(\kappa_0)} \Po & = \Po^\Sm \Go^{-1} \Ko^{(\kappa_0)} \Po^\Sm\\
    & + \paren{\iota \omega c_0^{-1} d}^{-1} \Po^\Sm \Go^{-1} \Ko^{(\kappa_0)} \Po^{\Lambda\textrm{H}} \, (=: \Ro^{sl})\\
    & + \paren{\iota \omega c_0^{-1} d} \Po^{\Lambda\textrm{H}} \Go^{-1} \Ko^{(\kappa_0)} \Po^\Sm \\
    & + \Po^{\Lambda\textrm{H}} \Go^{-1} \Ko^{(\kappa_0)} \Po^{\Lambda\textrm{H}}.
\end{align*}
Under infinite-precision quadrature and a conforming quasi-Helmholtz decomposition, \eqref{eq:vanishing_M0}, together with the algebraic structure of $\Go^{-1}$ in \eqref{eq:Go_inv} and the projector orthogonality $\Po^\Sm \Po^{\Lambda\textrm{H}} = \zrb$, implies that
\[
    \Po^\Sm \Go^{-1} \Ko^{(\kappa_0)} \Po^{\Lambda\textrm{H}} = \OO(\omega^2).
\]
Multiplication by $\paren{\iota\omega c_0^{-1} d}^{-1}$ yields $\Ro^{sl} = \OO(\omega)$. 

In a direct assembly of $\Ko^{(\kappa_0)}$, the dynamic and static contributions are evaluated together. Quadrature and discretization errors may then leave a small residual of the theoretically vanishing static term in the projected block. Since this residual is subsequently multiplied by a factor $\OO(\omega^{-1})$, it may produce a spurious $\OO(\omega^{-1})$ contribution and thereby compromise the low-frequency conditioning of the regularized system.

To eliminate this source of numerical instability, the static contribution is extracted analytically before quadrature. More precisely, the matrix $\Ro^{sl}$ is assembled using $\Ko_{\textrm{ext}}^{(\kappa_0)}$ in place of $\Ko^{(\kappa_0)}$, where $\Ko_{\textrm{ext}}^{(\kappa_0)}$ denotes the mixed discretization of the static-extracted operator $K_{\textrm{ext}}^{(\kappa_0)}$. Its kernel is obtained by replacing the original Green function by its dynamic remainder
\begin{equation}
    \label{eq:grad_extkernel}
    \nabla_{\xb} \dfrac{e^{-\iota \kappa_0 R} - 1}{4\pi R} = -\dfrac{\xb - \yb}{4\pi R^3} \paren{\paren{1 + \iota z} \paren{e^{-\iota z} - 1} + \iota z},
\end{equation}
with $z = \kappa_0 R$. Here, $e^{-\iota z} - 1$ is evaluated using the numerically stable \texttt{expm1} function available in most mathematical software libraries. The $\OO(\omega^2)$ behavior is thus enforced directly at the kernel level, rather than relying on the numerical cancellation of $\OO(1)$ contributions. Accordingly, the block $\Ro^{sl}$ retains the expected $\OO(\omega)$ scaling as $\omega \to 0$. 

The remaining three blocks do not contain an $\OO(\omega^{-1})$ rescaling and are not susceptible to the same amplification. They can therefore be assembled using the matrix $\Ko^{(\kappa_0)}$. 

An analogous treatment is required for the right-hand side. Any constant polarization vector $\pb$ admits a scalar potential $\phi$ such that $\pb = \nabla \phi$. Hence, for any solenoidal function $\gb$,
\[
    \inprod{\gb, \nv \times \pb}_{\times, \Gm} = \inprod{\gb, \nv \times \nabla \phi}_{\times, \Gm} = - \int_\Gm \phi \nabla_\Gm \cdot \gb \ds = 0.
\]
For the plane-wave excitation considered here, this orthogonality relation, together with the Taylor expansion of the exponential factor and the structure of $\Go^{-1}$, implies that
\begin{equation}
    \label{eq:loop_rhs}
    \abs{\Po^\Sm \Go^{-1} \eo^{inc}} = \OO(\omega), \qq \abs{\Po^\Sm \Go^{-1} \ho^{inc}} = \OO(\omega).
\end{equation}

The rescaled excitation vector in \eqref{eq:stabilized_MFIE} can be decomposed as
\[
    \Po^{-1} \Go^{-1} \ho^{inc} = \Po^{\Lambda\textrm{H}} \Go^{-1} \ho^{inc} + \paren{\iota \omega c_0^{-1} d}^{-1} \Po^\Sm \Go^{-1} \ho^{inc}.
\]
By \eqref{eq:loop_rhs}, the second component above, denoted by $\ro_s$, remains bounded as $\omega \to 0$. At the discrete level, however, the continuous orthogonality relation may not be satisfied exactly by the quasi-Helmholtz projection and the numerical quadrature. A residual static component in $\Po^\Sm\Go^{-1}\ho^{\mathrm{inc}}$ would therefore be amplified by the factor $\OO(\omega^{-1})$.

For this reason, $\ro_s$ is assembled using a static-extracted vector $\ho_{\textrm{ext}}^{inc}$, obtained by replacing $e^{-\iota \kappa_0 \wh{\kb} \cdot \xb}$ in \eqref{eq:plane_wave} with $e^{-\iota \kappa_0 \wh{\kb} \cdot \xb} - 1$, evaluated using the \texttt{expm1} function. The other component is assembled using the original vector $\ho^{inc}$. 


\subsection{M\"{u}ller Equation}
\label{subsec:Mueller_implementation}

The diagonal blocks and the excitation vector of the stabilized M\"uller formulation \eqref{eq:stabilized_Mueller} have the same structure as the regularized MFIE. They can therefore be assembled using the procedure described in Section~\ref{subsec:MFIE_implementation}. The off-diagonal blocks, however, require a different numerical treatment. Their low-frequency behavior relies on the cancellation of the hyper-singularities in the interior and exterior single-layer operators. This cancellation renders the complete off-diagonal blocks of order $\OO(\omega)$ and is fundamental to the low-frequency stability of both the original and regularized M\"uller formulations \cite{LC2026}.

At the continuous level, the interior and exterior operators are combined with weights chosen such that their static hyper-singular contributions cancel exactly. However, if the two operators are discretized independently, their common static components may be approximated differently. Quadrature errors can then leave a small residual after the two contributions are combined. Since the scalar-potential part of each single-layer operator contains the factor $(\iota \kappa_i)^{-1} (i = 0$ or $1)$, this residual may be amplified as the frequency decreases, resulting in a spurious loss of conditioning. 

To enforce the hyper-singularity cancellation at the discrete level, the common static contribution is removed analytically before quadrature. For the exterior operator, we introduce the static-extracted single-layer operator
\begin{align*}
    \paren{T^{(\kappa_0)}_{\textrm{ext}} \jb} (\xb) & := -\iota \kappa_0 \nv \times \int_{\Gm} \dfrac{e^{-\iota \kappa_0 R}}{4\pi R} \jb(\yb) \ds_{\yb} \\
    & \q\,\, + \dfrac{1}{\iota \kappa_0} \nv \times \nabla \int_{\Gm} \dfrac{e^{-\iota \kappa_0 R} - 1}{4\pi R} \nabla_\Gm \cdot \jb(\yb) \ds_{\yb}.
\end{align*}
Only the scalar-potential contribution is modified. The vector-potential contribution is left unchanged.

The gradient of the extracted kernel, evaluated using \eqref{eq:grad_extkernel}, satisfies
\[
    \nabla_{\xb} \dfrac{e^{-\iota \kappa_0 R}-1}{4 \pi R} = \OO(\kappa_0^2).
\]
Consequently, multiplication by $(\iota \kappa_0)^{-1}$ yields a scalar-potential contribution of order $\OO(\omega)$. The vector-potential term is also of order $\OO(\omega)$. The required $\OO(\omega)$ behavior is thus imposed directly at the kernel level and no longer depends on the numerical cancellation of separately assembled static hyper-singular terms.

The interior operator $T^{(\kappa_1)}$ is treated analogously by replacing it with $T_{\textrm{ext}}^{(\kappa_1)}$. The off-diagonal blocks are then assembled from the appropriately weighted combination of the two static-extracted operators. Because the common static hyper-singular contribution is removed analytically before quadrature, no further modification of the individual quasi-Helmholtz sub-blocks is required.

\subsection{MT-M\"{u}ller Equation}

In the stabilized MT-M\"uller formulation \eqref{eq:stablized_MT_Mueller}, the diagonal subdomain blocks $\MMs_{kk}$ and the excitation vector can be assembled using the procedures described in Sections~\ref{subsec:MFIE_implementation} and~\ref{subsec:Mueller_implementation}. Each off-diagonal subdomain block $\MMs_{jk}$, with $j \neq k$, is itself a $2 \times 2$ operator matrix. Its two off-diagonal entries can be assembled following Section~\ref{subsec:Mueller_implementation}, whereas the two diagonal entries require a separate treatment because they contain geometric identity operators that couple basis functions supported on the different boundaries $\Gm_k$ and $\Gm_j$.

A direct extension of the static MFIE cancellation property \eqref{eq:vanishing_M0} to the MT setting gives
\begin{equation}
    \label{eq:vanishing_M0_MT}
    \inprod{\gb, \paren{\dfrac{1}{2}I - K^{(0)}_{jk}} \fb}_{\times, \Gm_j} = 0,
\end{equation}
when $\fb$ and $\gb$ are both local-loop functions, or one is a local-loop function and the other is harmonic. Here, $\fb$ is supported on $\Gm_k$, whereas $\gb$ is defined on $\Gm_j$.

After mixed discretization, multiplication by the inverse Gram matrix, and application of the quasi-Helmholtz projectors, \eqref{eq:vanishing_M0_MT}, together with the algebraic structure of $\Go_j^{-1}$, implies
\begin{equation}
    \label{eq:dis_vanishing_MT}
    \Po^\Sm_j \Go_j^{-1} \paren{\dfrac{1}{2} \Go_{jk} - \Ko^{(0)}_{jk}} \Po^{\Lambda\textrm{H}}_k = \zrb.
\end{equation}
Here, $\Go_{jk}$ and $\Ko_{jk}^{(0)}$ are the mixed discretized matrices of the geometric identity operator $I$ and the static operator $K_{jk}^{(0)}$, respectively, employing the trial space $\mathrm{RT}(\Gm_{h,k})$ and the test space $\mathrm{BC}(\Gm_{h,j})$. As in the MFIE case, the cancellation \eqref{eq:dis_vanishing_MT} concerns only one of the four quasi-Helmholtz sub-blocks of each diagonal entry within $\MMs_{jk}$, which is denoted by $\Ro^{sl}_{jk}$.

For $j = k$, the mixed Gram matrix reduces to $\Go_{kk} = \Go_k$, so that \eqref{eq:dis_vanishing_MT} becomes
\[
    \dfrac{1}{2} \Po^\Sm_k \Po^{\Lambda\textrm{H}}_k - \Po^\Sm_k \Go_{k}^{-1} \Ko^{(0)}_{kk} \Po^{\Lambda\textrm{H}}_k = \zrb.
\]
Since $\Po^\Sm_k \Po^{\Lambda\textrm{H}}_k = \zrb$ by construction, it follows that
\[
    \Po^\Sm_k \Go_{k}^{-1} \Ko^{(0)}_{kk} \Po^{\Lambda\textrm{H}}_k = \zrb.
\]
Thus, within the diagonal subdomain blocks $\MMs_{kk}$, the projected identity and static contributions vanish independently. The required low-frequency scaling can therefore be enforced by extracting the static contribution from the kernel of $K_{kk}^{(\kappa_i)} (i = 0$ or $k)$, exactly as in the regularized MFIE.

The situation differs when $j \neq k$ and $\Gm_j \cap \Gm_k \neq \emptyset$. In this case, neither $\Po^\Sm_j \Go_j^{-1} \Go_{jk} \Po^{\Lambda\textrm{H}}_k$ nor $\Po^\Sm_j \Go_j^{-1} \Ko^{(0)}_{jk} \Po^{\Lambda\textrm{H}}_k$ vanishes individually. Therefore, extracting the static part of $\Ko^{(\kappa_i)}_{jk}$ alone would destroy this cancellation, thereby modifying the discrete operator and producing an incorrect solution.

Accordingly, in the quasi-Helmholtz sub-block $\Ro_{jk}^{sl}$ of each diagonal entry of the off-diagonal subdomain block $\MMs_{jk}$, the geometric identity contribution and the static part of $\Ko_{jk}^{(\kappa_i)}$ must be removed together. As in Section~\ref{subsec:MFIE_implementation}, the latter is eliminated analytically by replacing the Green kernel with its static-extracted counterpart. This joint treatment preserves the cancellation in \eqref{eq:dis_vanishing_MT} and prevents quadrature errors in the static contributions from being amplified by the low-frequency rescaling.

\subsection{Complexity Analysis}

We now assess the computational cost introduced by the proposed stabilizations. The mixed discretization of the MFIE and M\"uller formulations requires assembly on the barycentrically refined mesh and therefore increases the assembly cost by a constant factor of six relative to a discretization on the original mesh, such as the classical non-conforming scheme. Since this factor is independent of the number of unknowns, it does not affect the asymptotic complexity of the assembly.

The stabilization additionally requires the construction and application of the quasi-Helmholtz projectors and the action of the inverse Gram matrices. As shown in \cite{ACB+2013}, the quasi-Helmholtz projectors can be constructed and applied in linear complexity. Moreover, the Gram matrices are sparse and uniformly well conditioned. Their inverse can be computed in linear complexity using an iterative sparse solver.

Consequently, all additional operations introduced by the stabilization have at most linear complexity and are asymptotically dominated by the assembly and application of the boundary integral operators. The proposed regularizations therefore preserve the asymptotic complexity of the underlying MFIE and M\"uller formulations. In particular, when combined with a fast integral equation method, such as the fast multipole method or an appropriate hierarchical compression scheme, the stabilized formulations retain the quasi-linear computational complexity of the corresponding original formulations.

\section{Numerical Results}
\label{sec:results}


\subsection{Experimental Setup}

All numerical experiments are performed with free space as the background medium. The incident field is a plane wave with polarization $\pb = (1, 0, 0)^\transpose$ and direction $\wh{\kb} = (0, 0, 1)^\transpose$. 

\begin{figure}[!t]
    \centering
    \includegraphics[trim={0.1cm 0.6cm 0.2cm 0.6cm}, clip, width=\linewidth]{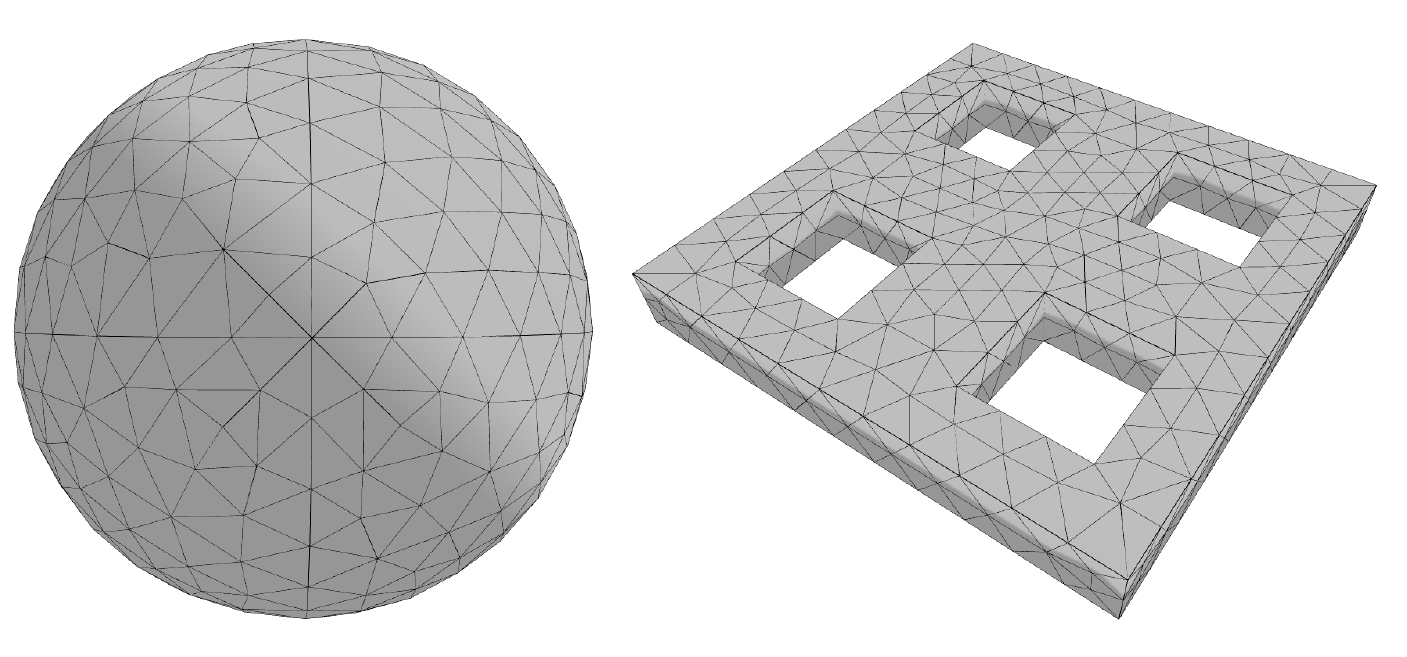}
    \caption{Scatterer geometries used in the PEC and homogeneous dielectric experiments. \textit{Left}: Unit sphere. \textit{Right}: Square torus with four holes.}
    \label{fig:mesh}
\end{figure}

The geometries employed in the PEC and homogeneous dielectric scattering experiments are shown in Fig.~\ref{fig:mesh}. The first scatterer is a sphere of radius $1 \, \mathrm{m}$. The second is a square torus with overall dimensions $1.6 \, \mathrm{m} \times 1.6 \, \mathrm{m} \times 0.2 \, \mathrm{m}$ and four square holes, each of dimensions $0.4 \, \mathrm{m} \times 0.4 \, \mathrm{m} \times 0.2 \, \mathrm{m}$. The sphere is smooth and simply connected, whereas the square torus is non-smooth and multiply connected with genus $g = 4$. For the far-field experiments, the sphere and the square torus are discretized with mesh sizes $h = 0.2 \, \mathrm{m}$ and $h = 0.15 \, \mathrm{m}$, respectively.

In the homogeneous dielectric experiments, the sphere and square torus are filled with materials of permittivity $\epsilon_1 = 3\epsilon_0$ and $\epsilon_1 = 100\epsilon_0$, respectively. The latter corresponds to a high-contrast configuration. In both cases, the permeability is set to $\mu_1 = \mu_0$.

\begin{figure}[!t]
    \centering
    \includegraphics[trim={0.1 0.6cm 0.2cm 0.6cm}, clip, width=\linewidth]{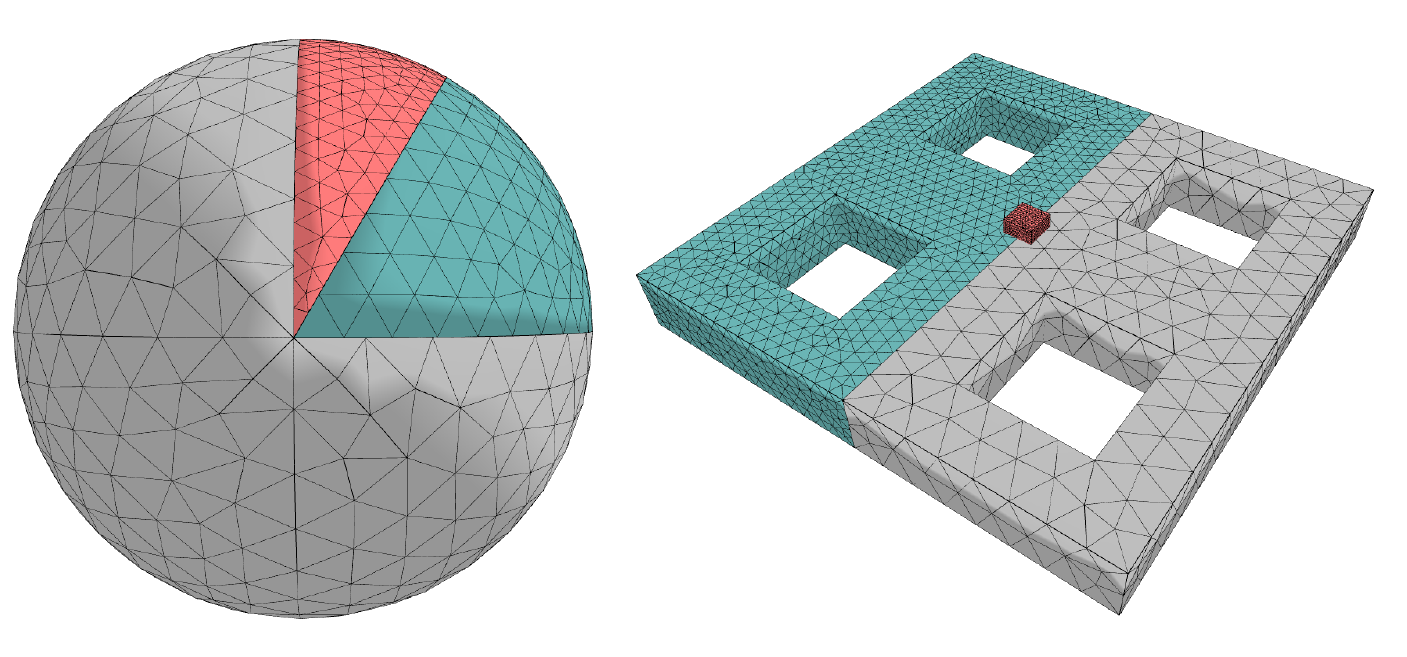}
    \caption{Composite dielectric scatterers discretized using non-conformal surface meshes. \textit{Left}: Unit sphere partitioned into three subdomains. \textit{Right}: Square torus partitioned into two subdomains, with a cuboidal boss attached to its upper surface.}
    \label{fig:nonconformalmesh}
\end{figure}

The geometries used in the composite dielectric experiments are depicted in Fig.~\ref{fig:nonconformalmesh}. The first composite scatterer is a unit sphere partitioned into three subdomains. The first subdomain occupies three quadrants of the sphere. The remaining quadrant is divided such that the second subdomain occupies two thirds of it and the third subdomain occupies the remaining third. The three subdomains are characterized by $(\epsilon_1, \epsilon_2, \epsilon_3) = (7 \epsilon_0, 14\epsilon_0, 42\epsilon_0)$ and $\mu_1 = \mu_2 = \mu_3 = \mu_0$. This configuration reproduces the experiment in \cite{LC2026}.

The second composite scatterer is based on the square torus introduced above. The torus is partitioned into two subdomains, and a boss of dimensions $0.1 \, \mathrm{m} \times 0.1 \, \mathrm{m} \times 0.05 \, \mathrm{m}$ is attached to its upper surface as a third subdomain. The three subdomains are assigned the material parameters $(\epsilon_1, \epsilon_2, \epsilon_3) = (2\epsilon_0, 10 \epsilon_0, 200 \epsilon_0)$ and $\mu_1 = \mu_2 = \mu_3 = \mu_0$. This configuration constitutes a multiscale, high-contrast test case. Both composite scatterers contain multiple junction edges, making them challenging for the treatment of coupled material interfaces. 

For the far-field experiments involving the composite sphere, the boundaries of the three subdomains are meshed independently. The mesh size in each subdomain boundary is chosen as $ h_k = 0.08 \lambda_k$, where $\lambda_k$ denotes the wavelength in the corresponding subdomain $\Om_k$ at $f = 50 \, \mathrm{MHz}$. For the composite square torus, the subdomain boundaries are also independently discretized, using $h_k = 0.05 \lambda_k$ at $f = 100 \, \mathrm{MHz}$.

The independently generated surface meshes are generally non-conformal across material interfaces. The interaction matrices are then evaluated using the quadrature strategy proposed in \cite{MC2026}, which constructs an on-the-fly local refinement for each interacting pair of triangles.

\subsection{Far-Field Evaluation}

\begin{figure*}[!t]
    \centering
    \includegraphics[trim={0 0.25cm 0 0.25cm}, clip, width=\linewidth]{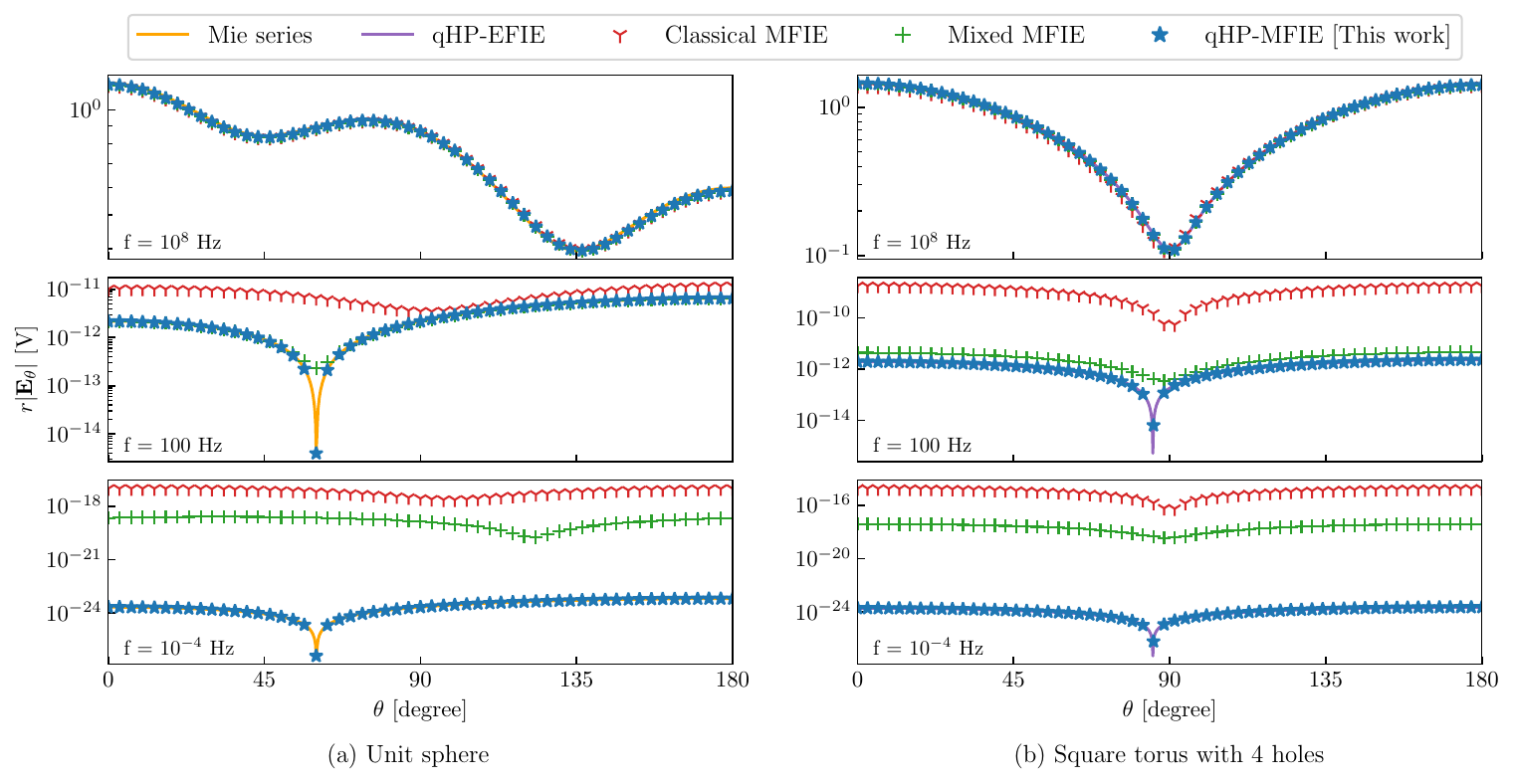}
    \caption{Far-field patterns for scattering from PEC objects at different frequencies, computed using various formulations. Results for the unit sphere are compared with the analytical Mie-series solution, whereas those for the square torus are compared with the qHP-EFIE solution.}
    \label{fig:MFIE_farfield}
\end{figure*}

We first consider electromagnetic scattering from the PEC sphere and square torus. Fig.~\ref{fig:MFIE_farfield} presents the far-field patterns computed using the classical non-conforming MFIE \eqref{eq:classical_MFIE}, the mixed MFIE \eqref{eq:dis_MFIE}, and the proposed stabilized MFIE \eqref{eq:stabilized_MFIE} (hereafter referred to as the qHP-MFIE). For the sphere, the numerical results are compared with the analytical Mie-series solution. For the square torus, the quasi-Helmholtz projected (qHP) EFIE proposed in \cite{ACB+2013} is used as the reference formulation. The frequencies range from the full-wave regime at $f = 100 \, \mathrm{MHz}$ to the low-frequency regime at $f = 100 \, \mathrm{Hz}$ and the extremely low-frequency regime at $f = 10^{-4} \, \mathrm{Hz}$. At $f = 100 \, \mathrm{MHz}$, all three MFIE formulations agree closely with the corresponding reference solutions. As the frequency decreases, however, the classical non-conforming MFIE progressively loses accuracy, followed at lower frequencies by the mixed MFIE. In contrast, the qHP-MFIE remains accurate over the entire frequency range considered. 

On the multiply connected square torus, the qHP-MFIE retains the topology-induced nullspace of the static MFIE and therefore becomes ill-conditioned as the frequency approaches zero. Nevertheless, Fig.~\ref{fig:MFIE_farfield}(b) shows no loss of far-field accuracy. In the static limit, the plane-wave excitation does not couple to the global-loop modes spanning the nullspaces of the exterior and interior MFIE operators \cite{CAO+2009a}. Although this nullspace is present at the continuous level, it is not preserved exactly by the proposed qHP-MFIE. Quadrature errors in the singular operator kernels perturb the zero singular values. By contrast, the right-hand side involves a smooth incident field and can be evaluated much more accurately. Its spurious projections are therefore much smaller than the perturbed singular values. Consequently, the associated components remain negligible, yielding a bounded current and an accurate far field.

\begin{figure*}[!t]
    \centering
    \includegraphics[trim={0 0.25cm 0 0.25cm}, clip, width=\linewidth]{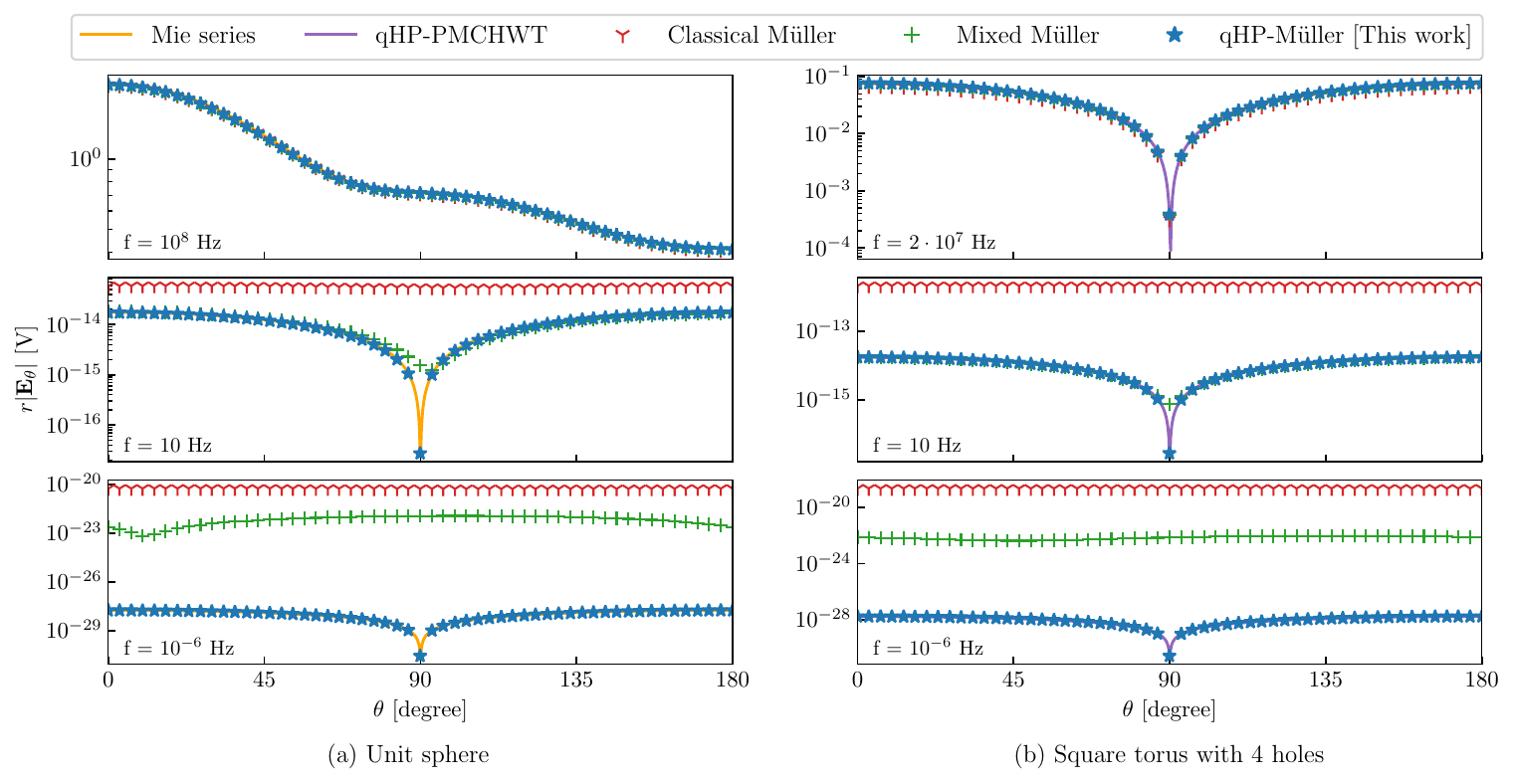}
    \caption{Far-field patterns for scattering from homogeneous dielectric objects at different frequencies, computed using various formulations. Results for the unit sphere are compared with the analytical Mie-series solution, whereas those for the square torus are compared with the qHP-PMCHWT solution.}
    \label{fig:Mueller_farfield}
\end{figure*}

We next examine scattering from homogeneous dielectric objects. Fig.~\ref{fig:Mueller_farfield} compares the far-field patterns obtained with the classical M\"uller formulation \eqref{eq:classical_Mueller}, the mixed M\"uller formulation \eqref{eq:dis_Mueller_1}, and the proposed stabilized formulation \eqref{eq:stabilized_Mueller} (denoted by qHP-M\"uller). As in the PEC experiments, the three formulations produce visually comparable and accurate results in the full-wave regime. At lower frequencies, the classical formulation loses accuracy because of the incorrect asymptotic scaling of its discretized system. As the frequency decreases further, cancellation errors also degrade the far-field predictions obtained with the mixed M\"uller formulation. The qHP-M\"uller formulation suppresses both sources of low-frequency inaccuracy and remains accurate over the entire frequency range considered. 
For the dielectric sphere, the analytical Mie-series solution is used as the reference. For the square torus, the results are compared with those obtained using the qHP-PMCHWT formulation proposed in \cite{BMC+2017}.

\begin{figure*}[!t]
    \centering
    \includegraphics[trim={0 0.25cm 0 0.25cm}, clip, width=\linewidth]{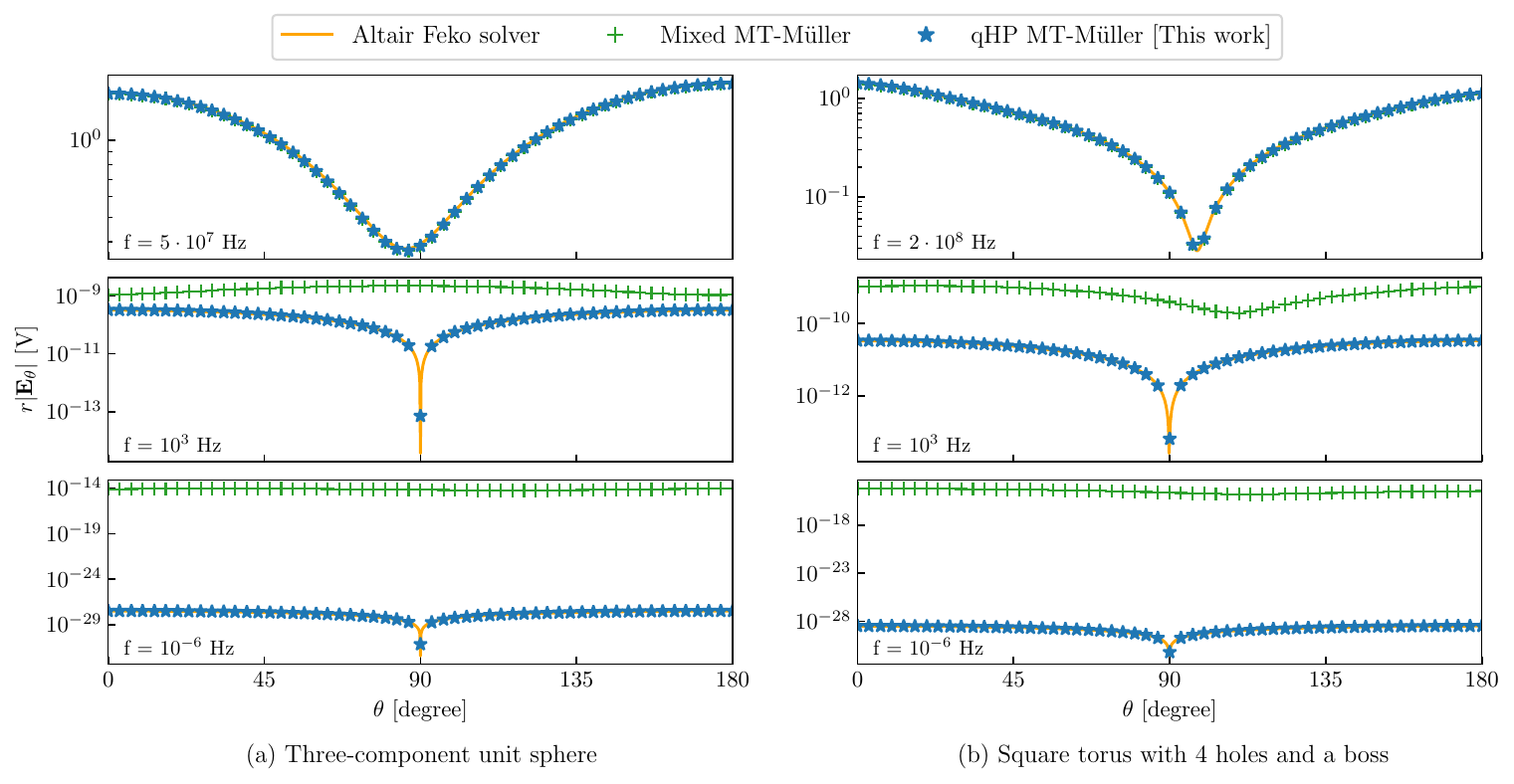}
    \caption{Far-field patterns for scattering from composite dielectric objects at different frequencies, computed using two MT-M\"uller formulations. Results are compared with reference solutions obtained using the electromagnetic solver Altair Feko 2026.}
    \label{fig:MT_Mueller_farfield}
\end{figure*}

Finally, we consider scattering from the composite dielectric objects shown in Fig.~\ref{fig:nonconformalmesh}. Fig.~\ref{fig:MT_Mueller_farfield} compares the mixed MT-M\"uller formulation \eqref{eq:dis_MT_Mueller} with the proposed stabilized qHP MT-M\"uller formulation \eqref{eq:stablized_MT_Mueller}. The mixed MT-M\"uller formulation begins to lose far-field accuracy at substantially higher frequencies than the mixed M\"uller formulation for a homogeneous dielectric object illustrated in Fig.~\ref{fig:Mueller_farfield}. In contrast, the qHP MT-M\"uller formulation remains accurate down to extremely low frequencies, even in the multiscale, high-contrast configuration. Reference solutions for these experiments are computed using the commercial electromagnetic solver Altair Feko 2026.

The frequency ranges displayed in Figs.~\ref{fig:MFIE_farfield}-\ref{fig:MT_Mueller_farfield} are selected to illustrate the onset of low-frequency cancellation errors in the unstabilized formulations and the effectiveness of the proposed regularizations. Numerical experiments conducted outside these ranges exhibit the same qualitative behavior.

\subsection{Conditioning}

This section investigates the conditioning of the proposed regularizations for the MFIE and the M\"uller formulations. In all experiments, the number of GMRES iterations required to reach the relative tolerance of $10^{-6}$ is measured.

\begin{figure}[!t]
    \centering
    \includegraphics[width=\linewidth]{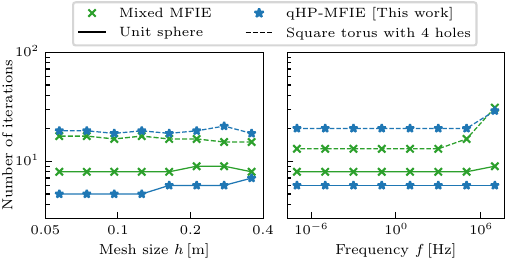}
    \caption{Number of GMRES iterations required by the mixed MFIE and the qHP-MFIE to reach the relative tolerance of $10^{-6}$ as functions of the mesh size $h$ and the frequency $f$.}
    \label{fig:MFIE_cond}
\end{figure}

Fig.~\ref{fig:MFIE_cond} shows the iteration counts of the mixed MFIE and the qHP-MFIE as functions of the mesh size $h$ and the frequency $f$. For the simply connected sphere, the iteration counts of both formulations remain nearly independent of mesh refinement and frequency in the low-frequency regime. The qHP-MFIE, however, consistently converges faster than the mixed MFIE. This behavior can be attributed to the structure of the stabilized system matrix in \eqref{eq:stabilized_MFIE}. After multiplication by the inverse Gram matrix and application of the rescaling operators, the system matrix takes the form of the identity plus a matrix representation of a compact operator, since $K^{(\kappa_0)}$ is compact on smooth surfaces. Its eigenvalues therefore cluster more tightly around unity, resulting in faster GMRES convergence. In contrast, the mixed MFIE matrix consists of the sparse, well-conditioned Gram matrix $\Go$ plus the discretization of $K^{(\kappa_0)}$. Although this system is also well conditioned, its spectrum is distributed over a broader region, leading to higher iteration counts.

For the multiply connected square torus, the condition numbers of both the mixed MFIE and the qHP-MFIE are expected to grow with decreasing frequency because of the topology-induced nullspace of the static MFIE operator \cite{CAO+2009a}. Even so, the experiments here show that the GMRES iteration counts remain nearly independent of both $h$ and $f$.
This may sound like a contradiction, but keep in mind that even though in the static limit the MFIE operator has an approximate nullspace of mesh-independent dimension $g$, its condition number away from this nullspace remains bounded.
%
The qHP-MFIE consistently requires more iterations than the mixed MFIE for this geometry because it partially preserves the cancellation properties of the static MFIE operators, yielding smaller singular values and hence slower convergence.

\begin{figure}[!t]
    \centering
    \includegraphics[width=\linewidth]{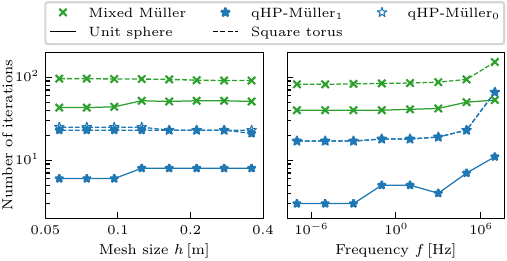}
    \caption{Number of GMRES iterations required by the mixed M\"uller formulation and the qHP-M\"uller formulation to reach the relative tolerance of $10^{-6}$ as functions of the mesh size $h$ and the frequency $f$. For the qHP-M\"uller formulation, the subscripts $0$ and $1$ indicate that the wave speeds $c_0$ and $c_1$, respectively, are used in the rescaling operator.}
    \label{fig:Mueller_cond}
\end{figure}

The corresponding results for the mixed M\"uller and qHP-M\"uller formulations are presented in Fig.~\ref{fig:Mueller_cond}. For both scatterer geometries, the iteration counts of the qHP-M\"uller formulation remain bounded under mesh refinement and in the low-frequency limit. Moreover, they are consistently lower than those of the mixed M\"uller formulation, indicating that the proposed regularization improves the conditioning of the mixed formulation, including for non-smooth geometries. For the high-contrast square torus, we additionally compare two variants of the qHP-M\"uller formulation in which the rescaling operators are constructed using either the exterior wave speed $c_0$ or the interior wave speed $c_1$. The corresponding iteration counts differ by at most one, demonstrating that this choice has a negligible effect on the convergence behavior.

\begin{figure}[!t]
    \centering
    \includegraphics[width=\linewidth]{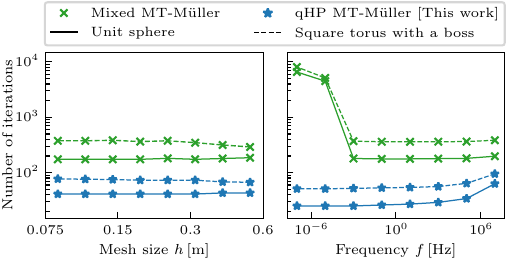}
    \caption{Number of GMRES iterations required by the mixed MT-M\"uller formulation and the qHP MT-M\"uller formulation to reach the relative tolerance of $10^{-6}$ as functions of the mesh size $h$ and the frequency $f$.}
    \label{fig:MT_Muller_cond}
\end{figure}

Finally, Fig.~\ref{fig:MT_Muller_cond} reports the iteration counts of the mixed MT-M\"uller and qHP MT-M\"uller formulations for the two composite dielectric scatterers. Over most of the investigated frequency range, the mixed MT-M\"uller formulation exhibits nearly constant iteration counts under mesh refinement and as the frequency decreases. At extremely low frequencies, however, the iteration count increases sharply, with the deterioration becoming apparent at approximately $f = 10^{-5} \, \mathrm{Hz}$ for both geometries. The cause of this abrupt loss of convergence has not yet been fully identified. In contrast, the qHP MT-M\"uller formulation maintains bounded and consistently lower iteration counts throughout the investigated ranges of mesh sizes and frequencies. This behavior is observed for both the composite sphere and the multiscale, high-contrast square torus, demonstrating that the proposed regularization remains effective in the presence of multiple interfaces with junctions, strong material contrasts, and multiple geometric scales.

\section{Conclusion}
\label{sec:conclusion}

We have developed numerical methods for suppressing low-frequency cancellation errors in the far-field evaluation with the MFIE and M\"uller formulations. The proposed methods employ rescaling operators constructed from quasi-Helmholtz projectors and appropriate frequency-dependent rescaling factors. The right rescaling introduces cancellation-free unknowns that preserve the Helmholtz components essential to far-field computation but otherwise susceptible to finite-precision cancellation at low frequencies. The left rescaling simultaneously compensates for the low-frequency ill-conditioning introduced by the right rescaling and balances the right-hand side. The resulting formulations preserve the favorable conditioning of the underlying mixed discretizations. The numerical results demonstrate that the proposed regularizations remain stable and accurate over a broad frequency range. Their effectiveness is maintained for smooth and non-smooth geometries, simply and multiply connected scatterers, homogeneous and composite dielectric objects, and configurations involving strong material contrasts, junction edges, and multiple geometric scales.

A natural direction for future work is the extension of the proposed framework to time-domain MFIE and M\"uller equations. Another promising direction is the development of a regularized second-kind formulation for composite scatterers containing both dielectric and PEC subdomains.


%




\ifCLASSOPTIONcaptionsoff
  \newpage
\fi

\IEEEtriggeratref{18}


\bibliographystyle{IEEEtran}
\bibliography{abrv_ref.bib}

%









\end{document}

%% file: lvc_style.tex
\newcommand{\R}{\mathbb R}

\newcommand{\G}{\mathbb G}
\newcommand{\M}{\mathbb M}

\newcommand{\Z}{\mathbb Z}

\newcommand{\Pbb}{\mathbb P}

\newcommand{\OO}{\mathcal O}

\newcommand{\ufrakbm}{\bm{\mathfrak u}}
\newcommand{\vfrakbm}{\bm{\mathfrak v}}

\DeclareMathOperator{\io}{i}
\DeclareMathOperator{\jo}{j}
\DeclareMathOperator{\mo}{m}
\DeclareMathOperator{\ho}{h}
\DeclareMathOperator{\eo}{e}
\DeclareMathOperator{\ro}{r}

\DeclareMathOperator{\Do}{D}
\DeclareMathOperator{\Io}{I}
\DeclareMathOperator{\Po}{P}

\DeclareMathOperator{\Go}{G}
\DeclareMathOperator{\Ko}{K}
\DeclareMathOperator{\To}{T}
\DeclareMathOperator{\Ro}{R}

\DeclareMathOperator{\diag}{diag}
\DeclareMathOperator{\range}{range}

\DeclareMathOperator{\vvs}{\mathbf v}

\DeclareMathOperator{\DDs}{\mathbf D}

\DeclareMathOperator{\PPo}{\mathbf P}
\DeclareMathOperator{\MMs}{\mathbf M}
\DeclareMathOperator{\GGs}{\mathbf G}

\newcommand{\paren}[1]{\left({#1}\right)}

\newcommand{\abs}[1]{\left\vert{#1}\right\vert}

\newcommand{\inprod}[1]{\left\langle{#1}\right\rangle}

\newcommand{\wt}[1]{\widetilde{#1}}
\newcommand{\wh}[1]{\widehat{#1}}
\newcommand{\ovl}[1]{\overline{#1}}

\newcommand{\transpose}{\mathsf{T}}

\newcommand{\bm}[1]{\boldsymbol{#1}}

\newcommand{\gm}{\gamma}
\newcommand{\Gm}{\Gamma}
\newcommand{\om}{\omega}
\newcommand{\Om}{\Omega}

\newcommand{\Sm}{\Sigma}

\newcommand{\sst}{\subset}

\newcommand{\Eb}{\bm{E}}

\newcommand{\Mb}{\bm{M}}

\newcommand{\xb}{\bm{x}}
\newcommand{\yb}{\bm{y}}

\newcommand{\hb}{\bm{h}}

\newcommand{\eb}{\bm{e}}

\newcommand{\ub}{\bm{u}}

\newcommand{\fb}{\bm{f}}
\newcommand{\gb}{\bm{g}}
\newcommand{\kb}{\bm{k}}
\newcommand{\jb}{\bm{j}}

\newcommand{\mb}{\bm{m}}

\newcommand{\pb}{\bm{p}}

\newcommand{\nv}{\textbf{n}}
\newcommand{\uv}{\textbf{u}}

\newcommand{\zrb}{\bm{0}}

\newcommand{\ds}{\,\mathrm{d}s}

\newcommand{\q}{\quad}
\newcommand{\qq}{\qquad}
\newcommand{\qqq}{\qquad\quad}
\newcommand{\qqqq}{\qquad\qquad}
\newcommand{\qqqqq}{\qquad\qquad\quad}

%% file: arXiv_qHP_MFIE_Mueller.bbl
\begin{thebibliography}{10}
\providecommand{\url}[1]{#1}
\csname url@samestyle\endcsname
\providecommand{\newblock}{\relax}
\providecommand{\bibinfo}[2]{#2}
\providecommand{\BIBentrySTDinterwordspacing}{\spaceskip=0pt\relax}
\providecommand{\BIBentryALTinterwordstretchfactor}{4}
\providecommand{\BIBentryALTinterwordspacing}{\spaceskip=\fontdimen2\font plus
\BIBentryALTinterwordstretchfactor\fontdimen3\font minus \fontdimen4\font\relax}
\providecommand{\BIBforeignlanguage}[2]{{%
\expandafter\ifx\csname l@#1\endcsname\relax
\typeout{** WARNING: IEEEtran.bst: No hyphenation pattern has been}%
\typeout{** loaded for the language `#1'. Using the pattern for}%
\typeout{** the default language instead.}%
\else
\language=\csname l@#1\endcsname
\fi
#2}}
\providecommand{\BIBdecl}{\relax}
\BIBdecl

\bibitem{HEA+2023}
B.~Hofmann, T.~F. Eibert, F.~P. Andriulli, and S.~B. Adrian, ``An excitation-aware and self-adaptive frequency normalization for low-frequency stabilized electric field integral equation formulations,'' \emph{{IEEE} Trans. Antennas Propag.}, vol.~71, no.~5, pp. 4301--4314, 2023.

\bibitem{ADC+2021}
S.~B. Adrian, A.~Dely, D.~Consoli, A.~Merlini, and F.~P. Andriulli, ``Electromagnetic integral equations: insights in conditioning and preconditioning,'' \emph{IEEE Open J. Antennas Propag.}, vol.~2, pp. 1143--1174, 2021.

\bibitem{ZCC+2003}
Y.~Zhang, T.~J. Cui, W.~C. Chew, and J.-S. Zhao, ``Magnetic field integral equation at very low frequencies,'' \emph{{IEEE} Trans. Antennas Propag.}, vol.~51, no.~8, pp. 1864--1871, 2003.

\bibitem{VGG+2013}
F.~Vico, Z.~Gimbutas, L.~Greengard, and M.~Ferrando-Bataller, ``Overcoming low-frequency breakdown of the magnetic field integral equation,'' \emph{IEEE Trans. Antennas Propag.}, vol.~61, no.~3, pp. 1285--1290, 2013.

\bibitem{LSC2018}
Q.~S. Liu, S.~Sun, and W.~C. Chew, ``A potential-based integral equation method for low-frequency electromagnetic problems,'' \emph{IEEE Trans. Antennas Propag.}, vol.~66, no.~3, pp. 1413--1426, 2018.

\bibitem{QC2010}
Z.-G. Qian and W.~C. Chew, ``Enhanced {A-EFIE} with perturbation method,'' \emph{{IEEE} Trans. Antennas Propag.}, vol.~58, no.~10, pp. 3256--3264, 2010.

\bibitem{ACB+2013}
F.~P. Andriulli, K.~Cools, I.~Bogaert, and E.~Michielssen, ``{On a well-conditioned electric field integral operator for multiply connected geometries},'' \emph{{IEEE} Trans. Antennas Propag.}, vol.~61, no.~4, pp. 2077--2087, 2013.

\bibitem{HEA+2023b}
B.~Hofmann, T.~F. Eibert, F.~P. Andriulli, and S.~B. Adrian, ``A low-frequency stable, excitation agnostic discretization of the right-hand side for the electric field integral equation on multiply-connected geometries,'' \emph{{IEEE} Trans. Antennas Propag.}, vol.~71, no.~12, pp. 9277--9288, 2023.

\bibitem{MBC+2020}
A.~Merlini, Y.~Beghein, K.~Cools, E.~Michielssen, and F.~P. Andriulli, ``{Magnetic and combined field integral equations based on the quasi-Helmholtz projectors},'' \emph{{IEEE} Trans. Antennas Propag.}, vol.~68, no.~5, pp. 3834--3846, 2020.

\bibitem{LCA+2024}
V.~C. Le, P.~Cordel, F.~P. Andriulli, and K.~Cools, ``{A stabilized time-domain combined field integral equation using the quasi-Helmholtz projectors},'' \emph{IEEE Trans. Antennas Propag.}, vol.~72, no.~7, pp. 5852--5864, 2024.

\bibitem{GAM+2017}
J.~E.~O. Guzman, S.~B. Adrian, R.~Mitharwal, Y.~Beghein, T.~F. Eibert, K.~Cools, and F.~P. Andriulli, ``{On the hierarchical preconditioning of the {PMCHWT} integral equation on simply and multiply connected geometries},'' \emph{{IEEE} Antennas Wirel. Propag. Lett.}, vol.~16, pp. 1044--1047, 2017.

\bibitem{BMC+2017}
Y.~Beghein, R.~Mitharwal, K.~Cools, and F.~P. Andriulli, ``{On a low-frequency and refinement stable PMCHWT integral equation leveraging the quasi-Helmholtz projectors},'' \emph{{IEEE} Trans. Antennas Propag.}, vol.~65, no.~10, pp. 5365--5375, 2017.

\bibitem{LMA+2026}
V.~C. Le, C.~M\"{u}nger, F.~P. Andriulli, and K.~Cools, ``{A stable, accurate, and well-conditioned time-domain PMCHWT formulation},'' \emph{IEEE Trans. Antennas Propag.}, 2026.

\bibitem{VFG+2016}
F.~Vico, M.~Ferrando, L.~Greengard, and Z.~Gimbutas, ``The decoupled potential integral equation for time‐harmonic electromagnetic scattering,'' \emph{Commun. Pure Appl. Math.}, vol.~69, no.~4, pp. 771--812, 2016.

\bibitem{CAY+2015}
J.~Cheng, R.~J. Adams, J.~C. Young, and M.~A. Khayat, ``{Augmented EFIE with normally constrained magnetic field and static charge extraction},'' \emph{IEEE Trans. Antennas Propag.}, vol.~63, no.~11, pp. 4952--4963, 2015.

\bibitem{CAO+2009a}
K.~Cools, F.~P. Andriulli, F.~Olyslager, and E.~Michielssen, ``{Nullspaces of {MFIE} and Calder{\'{o}}n preconditioned {EFIE} operators applied to toroidal surfaces},'' \emph{{IEEE} Trans. Antennas Propag.}, vol.~57, no.~10, pp. 3205--3215, 2009.

\bibitem{YTJ2005}
P.~Yl\"{a}-Oijala, M.~Taskinen, and S.~J\"{a}rvenp\"{a}\"{a}, ``Surface integral equation formulations for solving electromagnetic scattering problems with iterative methods,'' \emph{Radio Sci.}, vol.~40, no.~6, 2005.

\bibitem{YTJ2008}
------, ``Analysis of surface integral equations in electromagnetic scattering and radiation problems,'' \emph{Eng. Anal. Bound. Elem.}, vol.~32, no.~3, pp. 196--209, 2008.

\bibitem{UTR2011}
E.~Ubeda, J.~M. Tamay, and J.~M. Rius, ``Taylor-orthogonal basis functions for the discretization in method of moments of second kind integral equations in the scattering analysis of perfectly conducting or dielectric objects,'' \emph{Prog. Electromagn. Res.}, vol. 119, pp. 85--105, 2011.

\bibitem{KE2023}
J.~Kornprobst and T.~F. Eibert, ``Accuracy analysis of div-conforming hierarchical higher-order discretization schemes for the magnetic field integral equation,'' \emph{IEEE J. Multiscale Multiphysics Comput. Tech.}, vol.~8, pp. 261--268, 2023.

\bibitem{BC2007}
A.~Buffa and S.~H. Christiansen, ``A dual finite element complex on the barycentric refinement,'' \emph{Math. Comp.}, vol.~76, no. 260, pp. 1743--1770, 2007.

\bibitem{CAD+2011}
K.~Cools, F.~P. Andriulli, D.~De~Zutter, and E.~Michielssen, ``{Accurate and conforming mixed discretization of the {MFIE}},'' \emph{{IEEE} Antennas Wirel. Propag. Lett.}, vol.~10, pp. 528--531, 2011.

\bibitem{YJN2011}
S.~Yan, J.-M. Jin, and Z.~Nie, ``{Improving the accuracy of the second-kind Fredholm integral equations by using the Buffa-Christiansen functions},'' \emph{IEEE Trans. Antennas Propag.}, vol.~59, no.~4, pp. 1299--1310, 2011.

\bibitem{YJN2013}
------, ``Accuracy improvement of the second-kind integral equations for generally shaped objects,'' \emph{IEEE Trans. Antennas Propag.}, vol.~61, no.~2, pp. 788--797, 2013.

\bibitem{YKJ2016}
P.~Yla-Oijala, S.~P. Kiminki, and S.~Jarvenpaa, ``Conforming testing of electromagnetic surface-integral equations for penetrable objects,'' \emph{IEEE Trans. Antennas Propag.}, vol.~64, no.~6, pp. 2348--2357, 2016.

\bibitem{BCA+2014}
I.~Bogaert, K.~Cools, F.~P. Andriulli, and H.~Bagci, ``{Low-frequency scaling of the standard and mixed magnetic field and Müller integral equations},'' \emph{{IEEE} Trans. Antennas Propag.}, vol.~62, no.~2, pp. 822--831, 2014.

\bibitem{LC2026}
V.~C. Le and K.~Cools, ``{Multitrace M\"{u}ller boundary integral equation for electromagnetic scattering by composite objects},'' \emph{IEEE Trans. Antennas Propag.}, 2026.

\bibitem{RWG1982}
S.~Rao, D.~Wilton, and A.~Glisson, ``Electromagnetic scattering by surfaces of arbitrary shape,'' \emph{{IEEE} Trans. Antennas Propag.}, vol. AP-30, no.~3, pp. 409--418, 1982.

\bibitem{BCA+2011}
I.~Bogaert, K.~, Cools, F.~P. Andriulli, and D.~De~Zutter, ``Low frequency scaling of the mixed {MFIE} for scatterers with a non-simply connected surface,'' in \emph{Proc. Int. Conf. Electromagn. Adv. Appl. (ICEAA)}, 2011, pp. 951--954.

\bibitem{BCA2015b}
Y.~Beghein, K.~Cools, and F.~P. Andriulli, ``{A DC-stable, well-balanced, Calder{\'{o}}n preconditioned time domain electric field integral equation},'' \emph{{IEEE} Trans. Antennas Propag.}, vol.~63, no.~12, pp. 5650--5660, 2015.

\bibitem{Muller1969}
C.~M\"{u}ller, \emph{Foundations of the mathematical theory of electromagnetic waves}.\hskip 1em plus 0.5em minus 0.4em\relax Springer Berlin Heidelberg, 1969.

\bibitem{YT2005}
P.~Yla-Oijala and M.~Taskinen, ``{Well-conditioned M\"{u}ller formulation for electromagnetic scattering by dielectric objects},'' \emph{IEEE Trans. Antennas Propag.}, vol.~53, no.~10, pp. 3316--3323, 2005.

\bibitem{CHS2017}
X.~Claeys, R.~Hiptmair, and E.~Spindler, ``Second-kind boundary integral equations for electromagnetic scattering at composite objects,'' \emph{Comput. Math. Appl.}, vol.~74, no.~11, pp. 2650--2670, 2017.

\bibitem{Cools2026}
K.~Cools, ``{Fast converging single-trace quasi-local PMCHWT equation for the modeling of composite systems},'' \emph{IEEE Trans. Antennas Propag.}, vol.~74, no.~6, pp. 5697--5708, 2026.

\bibitem{CHJ2013}
X.~Claeys, R.~Hiptmair, and C.~Jerez-Hanckes, ``Multitrace boundary integral equations,'' in \emph{Direct and inverse problems in wave propagation and applications}, ser. Radon Ser. Comput. Appl. Math.\hskip 1em plus 0.5em minus 0.4em\relax Berlin: De Gruyter, 2013, vol.~14, pp. 51--100.

\bibitem{LBG+2022}
S.~Lasisi, T.~M. Benson, G.~Gradoni, M.~Greenaway, and K.~Cools, ``A fast converging resonance-free global multi-trace method for scattering by partially coated composite structures,'' \emph{IEEE Trans. Antennas Propag.}, vol.~70, no.~10, pp. 9534--9543, 2022.

\bibitem{MC2026}
C.~M\"{u}nger and K.~Cools, ``Global multi-trace, single source integral equation for the scattering by composite objects,'' \emph{IEEE Trans. Antennas Propag.}, 2026.

\end{thebibliography}
